\documentclass[11pt]{article}
\input epsf.tex
\usepackage{amssymb,latexsym,times}
\catcode`\@=11
\renewcommand\subsection{\@startsection{subsection}{2}{\z@}%
                                     {-3.25ex\@plus -1ex \@minus -.2ex}%
                                     {-0.01 mm}
                                     {\normalfont\large\bfseries}}

\catcode`\@=11
\renewcommand\subsubsection{\@startsection{subsubsection}{2}{\z@}%
                                     {-3.25ex\@plus -1ex \@minus -.2ex}%
                                     {-0.01 mm}
                                     {\normalfont\bfseries}}

\newtheorem{example}{Example}
\newtheorem{theorem}[example]{Theorem}
\newtheorem{corollary}[example]{Corollary}
\newtheorem{definition}[example]{Definition}
\newtheorem{proposition}[example]{Proposition}
\newtheorem{lemma}[example]{Lemma}
\newtheorem{conjecture}[example]{Conjecture}

\def\resp{{\em resp.$\ $}}

\def\proof{\medskip\noindent {\it Proof --- \ }}

\def\cqfd{\hfill $\Box$ \bigskip}
\def\adots{\mathinner{\mkern2mu\raise1pt\hbox{.}
\mkern3mu\raise4pt\hbox{.}\mkern1mu\raise7pt\hbox{.}}}
\def\<{\langle\,}
\def\>{\,\rangle}

\def\End{{\rm End}}
\def\Id{{\rm Id}}
\def\shuf{\sqcup\!\sqcup}
\def\ie{{\em i.e. }}
\def\eg{{\em e.g. }}

\def\SG{\mathfrak S}

\def\HHH{{\cal H}}

\def\l{\lambda}
\def\a{\alpha}
\def\de{\delta}
\def\b{\beta}
\def\ga{\gamma}
\def\N{{\mathbb N}}
\def\Z{{\mathbb Z}}
\def\C{{\mathbb C}}
\def\R{{\mathbb R}}
\def\Q{{\mathbb Q}}

\def\F{{\cal F}}

\def\B{{\bf B}}
\def\Bs{{\bf B}^*}
\def\U{{\cal U}}

\def\LL{{\mathfrak L}}

\def\SS{{\rm Supp}}

\def\ch{{\rm ch\, }}

\def\g{\mathfrak g}
\def\h{\mathfrak h}

\def\Sl{\mathfrak{sl}}

\def\nn{\mathfrak n}

\def\slchap{\widehat{\mathfrak{sl}}}

\def\A{{\cal A}}

\def\G{{\cal G}}
\def\LL{{\cal L}}
\def\L{\Lambda}

\def\<{\langle}
\def\>{\rangle}

\def\CC{{\cal C}}
\def\m{\mu}

\def\M{{\cal M}}
\def\RR{{\cal R}}

\def\le{\leqslant}
\def\ge{\geqslant}

\def\si{\sigma}

\def\De{\Delta}

\def\ga{\gamma}
\def\shuf{*}
\def\barshuf{\,\overline{*}\,}

\def\GL{{\cal GL}}
\def\D{\partial}
\def\k{\kappa}
\def\Ind{{\rm Ind}}
\def\mm{{\bf m}}
\def\ep{{\bf e'}}
\def\Ad{{\rm Ad}}

\def\shuff#1#2{\mathbin{
\hbox{\vbox{ \hbox{\vrule \hskip#2 \vrule height#1 width 0pt}%
\hrule}%
\vbox{ \hbox{\vrule \hskip#2 \vrule height#1 width 0pt
\vrule}%
\hrule}%
}}}

\def\SHUF{{\mathchoice{\shuff{7pt}{3.5pt}}%
{\shuff{6pt}{3pt}}%
{\shuff{4pt}{2pt}}%
{\shuff{3pt}{1.5pt}}}}%

\def\shuffle{\,\SHUF\,\,}

\title{\bf Dual canonical bases, quantum shuffles and $q$-characters} 

\author{Bernard {\sc Leclerc}}

\date{}

\begin{document}
\maketitle


\begin{abstract}\noindent
Rosso and Green have shown how to embed the positive part
$U_q(\nn)$ of a quantum enveloping algebra $U_q(\g)$ 
in a quantum shuffle algebra.
In this paper we study some properties of the image of the
dual canonical basis $\Bs$ of $U_q(\nn)$ under this embedding $\Phi$.
This is motivated by the fact that when $\g$ is of type $A_r$,
the elements of $\Phi(\Bs)$ are $q$-analogues of irreducible
characters of the affine Iwahori-Hecke algebras attached to the
groups $GL(m)$ over a $p$-adic field.
\end{abstract}

\section{Introduction}
In \cite{LLT1,LLT2,LT1} some close relationships were observed
between the representation theory of type $A$ Hecke algebras and
quantized Schur algebras on one side, and the canonical bases of certain
quantum groups on the other side.
Since then these connections have been studied by several authors
and other similar correspondences have been discovered 
\cite{A,LT2,VV,Gr,BK1,LNT,B,B2}.

Roughly speaking the principle is the following: the basis of simple 
modules in the Grothen\-dieck group of some appropriate category of 
representations of a certain algebra {\bf A} (\eg a Hecke algebra of type 
$A$ or $B$, an Ariki-Koike algebra,
a $q$-Schur algebra, $\ldots$) can be identified with the 
specialization at $q=1$ of the dual canonical basis of (a representation of) 
a certain quantum enveloping algebra {\bf U}. 

Let $\g$ be a complex simple Lie algebra and $\nn$ its maximal
nilpotent subalgebra.
In this paper we show that even when there is no apparent connection
with the Grothendieck group of some category, the dual canonical
basis $\Bs$ of $U_q(\nn)$ exhibits some features resembling certain classical
properties of the irreducible characters in Lie theory.

In order to observe this, one has to study $\Bs$ in the
particular realization of $U_q(\nn)$ discovered by Rosso \cite{R1,R2}
and Green \cite{G} in terms of quantum shuffles.
When $\g$ is of type $A_r$, the embedding $\Phi$ of $U_q(\nn)$
in the quantum shuffle algebra $\F$ is precisely a $q$-analogue
of the map $[M] \mapsto \ch M$ from the Zelevinsky ring
of a category of representations of the affine Hecke algebras
of type $GL(m)$ to the corresponding character ring. 
Indeed, as shown in \cite{Gr,GV}, the
multiplication of characters coming from parabolic induction
of Hecke modules is given by the (classical) shuffle product. 
Moreover, the above principle states in this case that 
the basis $\Bs$ is a $q$-analogue of the basis of the Zelevinsky ring consisting
of the classes of simple modules, hence $\Phi(\Bs)$ may be regarded as
a $q$-analogue of the set of irreducible characters
(see section~\ref{SECT4} below for more details).  

This motivates our investigation of $\Phi(\Bs)$ for general $\g$.
We shall be concerned with three main properties.
First, it is easy to describe explicitly the image $\Phi(U_q(\nn))$ of the
embedding of $U_q(\nn)$ in the quantum shuffle algebra
(Theorem~\ref{th2}).
We think of this result as an analogue of the classical
fact that the characters of the (virtual) integrable $\g$-modules
are the polynomials invariant under the action of the Weyl group.
Secondly, we show that the elements of $\Phi(\Bs)$ are parametrized
by their maximal word for the lexicographic order
(Theorem~\ref{th7}).
This is similar to the parametrization of irreducible 
integrable $\g$-modules by their highest weight.
Thirdly, it follows from Lusztig's geometric construction
of the canonical bases that for $\g$ of simply laced type,
the coefficients of the elements of $\Phi(\Bs)$ belong
to $\N[q,q^{-1}]$ (Theorem~\ref{positive}). 
We conjecture that this positivity property is also true in the 
non-simply laced case.
This is analogous to the fact that the character of a $\g$-module
is a positive sum of weights.

As an application, we describe in section~\ref{algo1} 
an algorithm for calculating the basis $\Bs$, which allowed us
to discover examples of imaginary vectors of $\Bs$ 
for $\g$ of type $A_5,B_3,C_3,D_4,G_2$, thus disproving a 
conjecture of Berenstein and Zelevinsky for all types except
$A_n\ (n\le 4)$ and $B_2$ \cite{Le}. 

There are some formal similarities between our results and
the theory of $q$-characters for finite-dimensional representations
of quantum affine algebras developed by Frenkel-Reshetikhin \cite{FR},
Frenkel-Mukhin \cite{FM} and Nakajima \cite{N}.
Actually, in type $A_r$, our $q$-characters for affine Hecke
algebras are interpreted geometrically in terms of the same
graded quiver varieties as those used by Nakajima for defining
the $(q,t)$-characters of $U_q(\slchap_n)$ \cite{L90}, so both
families of characters contain essentially the same information.
For other types though, there is no clear relationship
between $\Phi(\Bs)$ and the $q$-characters of $U_q(\widehat \g)$.

The paper is structured as follows.
In section~\ref{SECT1}, we review, following Rosso and Green,
the construction of the quantum shuffle embedding of $U_q(\nn)$.
Our presentation is based on the $q$-derivations $e'_i$ of
Kashiwara, which in the type $A$ case have the natural interpretation
of $i$-restriction operators in terms of affine Hecke algebras.
Then we prove Theorem~\ref{th2}.
We also describe explicitly the embedding of the algebra of regular 
functions $\C[N]$ in the (classical) shuffle algebra
obtained by specializing at $q=1$ the embedding $\Phi$
(here $N$ stands for a maximal unipotent subgroup of a complex simple 
Lie group $G$ with Lie algebra $\g$).
Sections~\ref{SECT2} and \ref{SECT2bis} are devoted to certain monomial
bases and PBW-type bases, respectively, which 
play an essential role in the proofs of our results on $\Bs$.
These two sections are based on some beautiful theorems of
Lalonde and Ram \cite{LR} and Rosso~\cite{R3}.
In particular, Lalonde and Ram have defined for any root system
a set of Lyndon words in one-to-one correspondence with the positive
roots.
These so-called `good Lyndon words' and their nonincreasing products
label in a natural way certain monomial and Lyndon bases.  
For the convenience of the reader, we have included proofs of 
most of the statements of \cite{LR,R3} needed for our purposes.
The main new result (Theorem~\ref{theomax}) describes the maximal words
of the images under $\Phi$ of the elements of certain Lusztig's
PBW-type bases. 
It is obtained by relating Lusztig's PBW-bases to Rosso's Lyndon bases.
In section~\ref{SECT3} we derive the above-mentioned properties
of $\Phi(\Bs)$ and we present an algorithm to compute it.
Section~\ref{SECT4} discusses the case of $\g$ of type $A_r$
and its relationship with the representation theory of affine
Hecke algebras, while section~\ref{SECT5} presents a 
conjectural analogue of this relationship for type $B_r$ and
the affine Hecke-Clifford superalgebras of Jones and Nazarov \cite{JN}
whose representation theory was recently studied by Brundan
and Kleshchev \cite{BK1}.
Finally, section~\ref{SECT6} describes a family of root vectors of
$\Bs$ for classical and simply-laced types.
More precisely, for the classical types $A_r, B_r, C_r, D_r$ we give a 
closed $q$-shuffle formula and for
the simply-laced types $A_r, D_r, E_r$ a simple combinatorial
description.
(For type $G_2$, the root vectors are calculated in \ref{EXG2}).   
This last section may serve to illustrate many statements of the paper.

\section{Embedding of $U_q(\nn)$ in a quantum shuffle algebra} 
\label{SECT1}

\subsection{}
Let $\g$ be a simple Lie algebra of rank $r$ over $\C$ and let $U_q(\g)$ be
the corresponding quantized enveloping algebra over $\Q(q)$ with Chevalley
generators $e_i, f_i  \ (i=1,\ldots ,r)$.
The Cartan matrix of $\g$ is denoted by $[a_{ij}]_{i,j=1,\ldots,r}$.
Let $\Delta$ be the root system of $\g$, $\Delta^+$ the subset of positive
roots, $Q$ the root lattice, $\Pi=\{\a_1,\ldots, \a_r\}$ the set of simple roots, 
$Q^+=\oplus_{i=1}^r \N\a_i$ the monoid generated by the simple roots, 
and $(\cdot\, , \cdot)$ a symmetric bilinear form on $Q$ such that 
\[
a_{ij} = {2(\a_i\,,\a_j)\over (\a_i\,,\a_i)} = {(\a_i\,,\a_j)\over
  d_i}, \qquad (1\le i,j \le r)
\]
where $d_i = (\a_i\,,\a_i)/2 \in \{1,2,3\}$. 

Let $U_q(\nn)$ be the subalgebra of $U_q(\g)$ generated by the
elements $e_i \ (i=1,\ldots ,r)$.
The defining relations of $U_q(\nn)$ are the so-called $q$-Serre
relations:
\begin{equation}
  \label{qSerre}
\sum_{k+l=1-a_{ij}} (-1)^k 
\left[
\matrix{1-a_{ij}\cr k}
\right]_i
e_i^k e_j e_i^l
=0, \qquad (1\le i \not = j \le r).
\end{equation}
Here we use the standard notation for $q$-integers and $q$-binomial
coefficients, namely,
\[
q_i=q^{d_i},\quad
[k]_i={q_i^k-q_i^{-k}\over q_i-q_i^{-1}}, \quad
\left[
\matrix{m\cr k}
\right]_i
={[m]_i[m-1]_i\cdots [m-k+1]_i\over[k]_i[k-1]_i\cdots[1]_i}.
\]
The algebra $U_q(\nn)$ is $Q^+$-graded by assigning to $e_i$ the degree $\a_i$.
We shall denote by $|u|$ the $Q^+$-degree of a homogeneous element
$u$ of $U_q(\nn)$.  

\subsection{}\label{ss1.2}
Kashiwara \cite{K} has introduced some $q$-derivations $e'_i  \ (i=1,\ldots ,r)$
of $U_q(\nn)$. 
These are the elements of $\End\, U_q(\nn)$ characterized by
\begin{equation}
  \label{eq:1.1}
e'_i(e_j)=\delta_{ij},\qquad
e'_i(uv)=e'_i(u)v+q^{-(\a_i,|u|)}ue'_i(v),
\end{equation}
for all homogeneous elements $u, v$ of $U_q(\nn)$.
It is known \cite{K} that 
\begin{equation}
  \label{eq:E'}
  e'_i(u)=0, \qquad (i=1,\ldots ,r)
\quad\Longleftrightarrow\quad
|u|=0.
\end{equation}
It is also known \cite{K} that these endomorphisms satisfy the $q$-Serre
relations, that is,
\begin{equation}
  \label{eq:serre}
\sum_{k+l=1-a_{ij}} (-1)^k 
\left[
\matrix{1-a_{ij}\cr k}
\right]_i
(e'_i)^k e'_j (e'_i)^l
=0, \qquad (1\le i \not = j \le r).
\end{equation}
Kashiwara \cite{K} proves that there is a unique nondegenerate
symmetric bilinear form
$(\cdot\,,\cdot)$ on $U_q(\nn)$ such that $(1,1)=1$ and 
\begin{equation}
  \label{eq:1.2}
(e'_i(u),v) = (u,e_iv), \qquad (u,v\in U_q(\nn),\ \ i=1,\ldots ,r),
\end{equation}
that is, $e'_i$ is the endomorphism adjoint to left multiplication
by $e_i$.

Note that Lusztig uses a slightly different scalar product $(\cdot\,,\cdot)_L$
satisfying
\begin{equation}
  \label{eq:1.22}
(e'_i(u),v)_L = {1\over 1-q^{2d_i}}\,(u,e_iv)_L, 
\qquad (u,v\in U_q(\nn),\ \ i=1,\ldots ,r),
\end{equation}
(see \cite{L}, 1.2.3, 1.2.13).
It is easy to see that if $u$ and $v$ are homogeneous elements
of $U_q(\nn)$ we have $(u,v)=(u,v)_L=0$ if $|u|\not = |v|$,
and 
\[
(u,v)_L=\prod_{i=1}^r {1\over (1-q^{2d_i})^{c_i}}\,(u,v)
\]
if $|u| = |v| = \sum_i c_i \a_i$.
It follows that if ${\cal B}$ is a basis of $U_q(\nn)$ consisting
of homogeneous vectors, then the adjoint bases of ${\cal B}$
with respect to $(\cdot\,,\cdot)$ and $(\cdot\,,\cdot)_L$
differ only by some normalization factors. 
In particular, ${\cal B}$ is orthogonal with respect to
$(\cdot\,,\cdot)$ if and only if it is orthogonal with respect to
$(\cdot\,,\cdot)_L$.
In this paper we shall use Kashiwara's form $(\cdot\,,\cdot)$.

\subsection{}
Let $\M$ (\resp $\F$) be the free monoid (\resp the free associative
algebra over $\Q(q)$) generated by the set of letters $I=\{w_1,\ldots ,w_r\}$.
We will use the notation 
$w[i_1,\ldots,i_k] := w_{i_1}\cdots w_{i_k}$.
The empty word is written $w[\,]$.
The length of a word $w\in\M$ is denoted by $\ell(w)$.
The algebra $\F$ is $Q^+$-graded by assigning to $w_i$ the degree
$\a_i$. The degree of a homogeneous element $f\in\F$ is denoted
by $|f|$.

In \cite{L}, Lusztig has endowed $\F$ with a twisted bialgebra
structure defined in terms of the bilinear form on $Q$.
He has shown that there exists a unique symmetric bilinear
form on $\F$ which adjoins the multiplication and the twisted
comultiplication.
Moreover the radical of this form coincides with the kernel
of the homomorphism $\F \longrightarrow U_q(\nn)$ mapping
$w_i$ to $e_i$, and the form it induces on $U_q(\nn)$ is nothing else
than the form (\ref{eq:1.22}) above.
Similarly $U_q(\nn)$ is endowed with a twisted bialgebra 
structure whose comultiplication is adjoint to the multiplication 
with respect to (\ref{eq:1.22}).
Hence by taking graded duals, we obtain a natural embedding
of vector spaces
\[
U_q(\nn) \cong U_q(\nn)^* \longrightarrow \F^* \cong \F
\]
in which the multiplication of $U_q(\nn)$ is sent to the
multiplication of $\F^*$ coming from Lusztig's comultiplication
on $\F$, and is a $q$-analogue of the shuffle product
as explained very clearly by Green \cite{Gr}.

Here we are going to indicate briefly how to recover this result 
by means of the $q$-derivations $e'_i$. 
An advantage of this approach is that it shows immediately
how this embedding specializes at $q=1$ to an embedding of
$\C[N]$ in the shuffle algebra, given explicitly in terms
of differential operators (see \ref{special1}).

\subsection{}
To $w=w[i_1,\ldots,i_k]$ we associate 
$\D_w:=e'_{i_1}\cdots e'_{i_k} \in \End\,U_q(\nn)$.
(For $w=w[\,]$ we set $\D_w = \Id_{U_q(\nn)}$.)
If $u$ is a homogeneous element of $U_q(\nn)$ and
$|w|=|u|$ then $\D_w(u)$ is of degree $0$, that is,
a scalar.
We define a $\Q(q)$-linear map 
$\Phi : U_q(\nn) \longrightarrow \F$
by setting
\begin{equation}
  \label{eq:1.3}
\Phi(u) = \sum_{w\in\M,\,|w|=|u|} \D_w(u) w  
\end{equation}
for a homogeneous element $u\in U_q(\nn)$.
It follows easily from (\ref{eq:E'}) that the map $\Phi$ is injective.

It may be helpful to think of (\ref{eq:1.3}) as a formal
Taylor expansion of $u$ 
(see below Proposition~\ref{varphi}).

\subsection{}
Define inductively a bilinear map $\shuf$ from $\F$ to $\F$ by setting,
for $a,b\in I$ and $w,x\in \M$
\begin{equation}
  \label{eq:1.4}
wa\shuf xb = (w\shuf xb)a +q^{-(|wa|,|b|)} (wa \shuf x)b,
\qquad
w[\,]\shuf x = x\shuf w[\,] = x.  
\end{equation}
Iterating (\ref{eq:1.4}) we get
\begin{equation}
  \label{eq:1.5}
  w[i_1,\ldots ,i_m]\shuf w[i_{m+1},\ldots ,i_{m+n}]
=
\sum_\sigma q^{-e(\sigma)} w[i_{\sigma(1)},\ldots ,i_{\sigma(m+n)}] 
\end{equation}
where 
the sum runs over the $\sigma \in \SG_{m+n}$ such that
$\sigma(1)<\cdots < \sigma(m)$ and $\sigma(m+1)<\cdots < \sigma(m+n)$,
and 
\begin{equation}
  \label{eq:1.6}
e(\sigma) = \sum_{k\le m<l;\ \sigma(k)<\sigma(l)}
(\a_{i_{\sigma(k)}}\,,\a_{i_{\sigma(l)}}).
\end{equation}
Thus, for $q=1$, $\shuf$ is the 
classical shuffle product $\shuffle$ in $\F$ \cite{Re},
and in particular it is associative and commutative.
The next proposition follows easily from the definitions
and its proof will be omitted.
\begin{proposition}\label{propo2}
The product $\shuf$ is associative, and for $w, x \in \M$
\[
w\shuf x = q^{-(|w|,|x|)} x\barshuf w
\]
where $\barshuf$ is the map obtained by replacing $q$ by $q^{-1}$
in the definition of $\shuf$. \cqfd
\end{proposition}

The following Lemma is a simple rank 2 computation.
\begin{lemma}\label{lemma3}
For $i\not = j$, 
\[
\sum_{k+l=1-a_{ij}} (-1)^k 
\left[
\matrix{1-a_{ij}\cr k}
\right]_i
w_i^{\shuf k} \shuf w_j \shuf w_i^{\shuf l}
=0,
\]
where we have put 
$w_i^{\shuf k} = w_i\shuf w_i \shuf \cdots \shuf w_i$ 
($k$ factors).
\cqfd
\end{lemma}

We introduce $\ep_i \in \End\,\F\ (i=1,\ldots,r)$ by setting
\begin{equation}
\ep_i(w[i_1,\ldots,i_k]) = \delta_{i,i_k} w[i_1,\ldots ,i_{k-1}],
\qquad \ep_i(w[\,])=0.
\end{equation}
\begin{lemma}\label{lemma4}
The endomorphisms $\ep_i$ satisfy
\[
\ep_i(w_j)=\delta_{ij},\qquad
\ep_i(x\shuf z)=\ep_i(x)\shuf z+q^{-(\a_i,|x|)}x\shuf \ep_i(z),
\]
for all homogeneous elements $x, z$ of $\F$.
\end{lemma}
\proof
Follows immediately from (\ref{eq:1.4}). \cqfd

\begin{theorem}[\cite{R1,R2,G}]\label{th1}
For $u,v \in U_q(\nn)$ we have $\Phi(uv)=\Phi(u)\shuf\Phi(v)$.
\end{theorem}
\proof
By Lemma~\ref{lemma3}, there exists a linear map 
$\Psi : U_q(\nn) \longrightarrow \F$ such that
\[
\Psi(e_i)=w_i\quad (i=1,\ldots, r),\qquad 
\Psi(uv) = \Psi(u)\shuf\Psi(v)\quad (u,v\in U_q(\nn)).   
\]
By Lemma~\ref{lemma4} this map satisfies:
$
\Psi e'_i = \ep_i \Psi, \ ( i=1,\ldots, r). 
$
Let $u\in U_q(\nn)$ be homogeneous and let $w=w[i_1,\ldots,i_k]\in\M$
be such that $|w|=|u|$.  
Let $\ga_w(u)$ be the coefficient of $w$ in $\Psi(u)$.
Then
$
\ga_w(u)=\ep_{i_1}\cdots \ep_{i_k}\Psi(u)
= \Psi e'_{i_1}\cdots e'_{i_k}(u) 
= \D_w(u).
$
Hence $\Psi(u) = \Phi(u)$, which proves the theorem.\cqfd

\subsection{}
By Theorem~\ref{th1}, the algebra $U_q(\nn)$ is isomorphic to the
subalgebra $\U$ of $(\F,\shuf)$ generated by the letters $w_i\in I$.
The next theorem gives a more explicit description of $\U$.

Let $i\not = j$ and $0\le k \le 1-a_{ij}$.
For $z,t\in\M$, we set 
$w(i,j,k;z,t) = zw_i^kw_jw_i^{1-a_{ij}-k}t$.
\begin{theorem}\label{th2}
The element $f=\sum_{w\in\M} \gamma(w) \, w$ of $\F$ belongs to
$\U$ if and only if
\begin{equation}\label{eqembed}
\sum_{k=0}^{1-a_{ij}} (-1)^k 
\left[
\matrix{1-a_{ij}\cr k}
\right]_i
\gamma(w(i,j,k;z,t))
=0
\end{equation}
for all $i\not = j$ and $z,t\in\M$. 
\end{theorem}
\proof
Let $K$ be the subspace of $\F$ defined by the system of linear
equations (\ref{eqembed}).
Let 
\[
f=\Phi(u)=\sum_{|w|=\nu}\gamma(w)\,w
\] 
for some $u\in U_q(\nn)$ of degree $\nu\in Q^+$.
Then for $w=w[i_1,\ldots ,i_k]$ of degree $\nu$ we have
\[
\gamma(w)=e'_{i_1}\cdots e'_{i_k}(u) 
= (e_{i_k}\cdots e_{i_1}\,,\,u)\,.
\]
Hence the fact that the elements $e_i$ satisfy the $q$-Serre relations 
(\ref{qSerre})
implies that $f\in K$. So $\U\subset K$.

Let $\F_\nu$ (\resp $K_\nu$, $U_q(\nn)_\nu$, $\U_\nu$)
be the homogeneous component of degree $\nu$ of $\F$
(\resp $K$, $U_q(\nn)$, $\U$).
Since (\ref{qSerre}) is a presentation of $U_q(\nn)$, we see that 
$\dim K_\nu = \dim U_q(\nn)_\nu$ for every $\nu\in Q^+$.
Moreover, $\Phi$ being injective, $\dim U_q(\nn)_\nu = \dim \U_\nu$,
hence $K=\U$.
\cqfd

\subsection{} \label{automorph}
The next proposition shows that some important automorphisms
of $U_q(\nn)$ can be seen as restrictions of certain simple 
linear maps defined over $\F$.
\begin{proposition}\label{proptau}
 {\rm (i)} Let $\tau$ be the $\Q(q)$-linear map from $\F$ to $\F$ 
such that 
\[
\tau(w[i_1,\ldots,i_k])=w[i_k,\ldots,i_1].
\]
Then $\tau(f \shuf g) = \tau(g) \shuf \tau(f)$
for all $f, g \in\F$.
Hence $\tau$ restricts to the $\Q(q)$-linear anti-automorphism
of $\U$ fixing the generators $w_i$.

\smallskip\noindent
{\rm (ii)} Let $f\mapsto \overline{f}$ be the $\Q$-linear map 
from $\F$ to $\F$ such that 
\[
\overline{q}=q^{-1},\qquad
\overline{w[i_1,\ldots,i_k]}=q^{-\sum_{1\le s<t\le k}(\a_{i_s},\,\a_{i_t})}\,
w[i_k,\ldots,i_1].
\]
Then $\overline{f \shuf g} = \overline{f} \shuf \overline{g}
\,$
for all $f,g$ in $\F$.
Hence $f\mapsto \overline{f}$ restricts to the $\Q$-linear automorphism
of $\U$ sending $q$ to $q^{-1}$ and fixing the generators $w_i$.

\smallskip\noindent
{\rm (iii)} Let $\si$ be the $\Q$-linear map 
from $\F$ to $\F$ such that 
\[
\si(q)=q^{-1},\qquad
\si(w[i_1,\ldots,i_k])=q^{-\sum_{1\le s<t\le k}(\a_{i_s},\,\a_{i_t})}\,
w[i_1,\ldots,i_k].
\]
Then $\si(f)=\overline{\tau(f)}$. Hence, $\si$ restricts to the $\Q$-linear 
anti-automorphism of $\U$ sending $q$ to $q^{-1}$ and fixing the 
generators $w_i$.
\end{proposition}
\proof It is enough to check (i) and (ii) when $f$ and $g$ are
two words. 
Then (i) follows immediately from (\ref{eq:1.5}).
To prove (ii) we may argue by induction on the length of the 
words. 
First note that 
\[
\overline{w[i_1,\ldots,i_k]}=q^{-(\a_{i_1},\,\a_{i_2}+\cdots
  +\a_{i_k})}\,
\overline{w[i_2,\ldots,i_k]}w_{i_1}.
\]
Assume by induction that (ii) is proved for every pair of words
whose sum of lengths is equal to $n$, and let 
$a,b\in I$ and $w,x\in \M$ with $\ell(aw)+\ell(bx)=n+1$.
Using (\ref{eq:1.5}), we have
\[
aw\shuf bx = q^{-(|a|,|bx|)}a(w\shuf bx) + b(aw \shuf x),
\]
hence
\begin{eqnarray*}
\overline{aw\shuf bx} &=& \overline{q^{-(|a|,|bx|)}a(w\shuf bx) + b(aw
  \shuf x)}\\
&=& q^{-(|a|,\,|w|)}\,(\overline{w\shuf bx})a
+ q^{-(|b|,\,|aw|+|x|)}\,(\overline{aw\shuf x})b \\
&=& q^{-(|a|,\,|w|)-(|b|,\,|x|)}\,(\overline{w}\shuf \overline{x}b)a
+ q^{-(|a|,\,|w|)-(|b|,\,|aw|+|x|)}\,(\overline{w}a\shuf \overline{x})b \\
&=& q^{-(|a|,\,|w|)-(|b|,\,|x|)}\,(\overline{w}a\shuf \overline{x}b)\\
&=& \overline{aw}\shuf \overline{bx}.
\end{eqnarray*}
Finally, (iii) follows from (i) and (ii).
\cqfd

For $\nu=\sum_{i=1}^r c_i\a_i\in Q^+$, define 
\begin{equation}
N(\nu)={1\over 2}\left((\nu\,,\nu)-\sum_{i=1}^r c_i (\a_i\,,\a_i)\right).
\end{equation}
Then, by Proposition~\ref{proptau}, 
$\sigma(w)=q^{-N(|w|)}\,w$ for all $w$,
from which the next lemma follows.
\begin{lemma}\label{symq}
Let $f=\sum_{w\in\M} \gamma_w(q)\, w\in\F$ be homogeneous. 
Then 
\[
\si(f)=q^{-N(|f|)} f\quad \Longleftrightarrow \quad 
\gamma_w(q^{-1})=\gamma_w(q), \quad (w\in\M).
\] 
\cqfd
\end{lemma} 

\subsection{}\label{special1}
We close this section by discussing the specialization 
of the $q$-shuffle embedding $\Phi$ at $q=1$.
This will be useful in sections~\ref{SECT4} and \ref{SECT5},
when we study characters of Hecke algebras.

\subsubsection{}\label{special2}
Let $\A=\Z[q,q^{-1}]$.
Following Lusztig, we introduce the $\A$-subalgebra $U_\A$
of $U_q(\nn)$ generated by the divided powers 
$e_i^{(k)}=e_i^k/[k]_i!\ (1\le i\le r,\ k\in \N)$.
We set 
\[
U^*_\A=\{u\in U_q(\nn) \ |\ (u\,,\,v)\in \A \mbox{ for all } v\in U_\A\}\,.
\]
For $w=w_{i_1}^{a_1}\cdots w_{i_k}^{a_k}\in\M$ with 
$i_j\not = i_{j+1}\ (1\le j\le k-1)$, we define 
$c_w=[a_1]_{i_1}!\cdots [a_k]_{i_k}!$
and we write 
$e_w=e_{i_1}\cdots e_{i_k}$. 
Thus $c_w^{-1}e_w$ is a product of divided powers.
Consider the free $\A$-module
$\F_\A = \bigoplus_{w\in\M} \A\,c_w w$ and 
define $\U_\A^*=\U\cap\F_\A$.
\begin{lemma}\label{Abasis}
We have $\U_\A^*=\Phi(U^*_\A)$.
\end{lemma}
\proof 
An element $u\in U_q(\nn)$ belongs to $U^*_\A$ if and only if
$(u\,,\,c_w^{-1}e_w)\in \A$ for all $w\in\M$, that is, if and
only if $\Phi(u)$ is an $\A$-linear combination of the elements
$c_w w$, that is, if and only if $\Phi(u)\in\F_\A$.
\cqfd

It is easy to see from (\ref{eq:1.5}) (\ref{eq:1.6}) that 
$\F_\A$ is in fact a subalgebra of $\F$.
It follows that $\U^*_\A$ is a subalgebra of $\U$, and by 
Lemma~\ref{Abasis} that $U^*_\A$ is a subalgebra of $U_q(\nn)$,
a well-known fact.
Define 
\[
\F_\C=\C\otimes_\A \F_\A, \qquad 
U^*_\C=\C\otimes_\A U^*_\A,\qquad
\U^*_\C=\C\otimes_\A \U^*_\A,
\]
where $\C$ is regarded as an $\A$-module via $q\mapsto 1$.
The natural maps 
$
\F_\A \rightarrow \F_\C
$
and
$
U^*_\A \rightarrow U^*_\C
$
will be called `specialization at $q=1$'.

The $\C$-linear map defined by 
$1\otimes c_w w \mapsto a_1!\cdots a_k!w$ is an algebra isomorphism
from $\F_\C$ endowed with the specialization of $\shuf$ at
$q=1$ to the classical $\C$-shuffle algebra over 
$\{w_1,\ldots ,w_r\}$, and from now on these two algebras
will be identified. 
The subalgebra $\U^*_\C$ of $\F_\C$ can be described explicitly
by specializing in the obvious way Theorem~\ref{th2} at $q=1$.
Note that $\U^*_\C$ is in general strictly bigger than the
subalgebra of $(\F_\C, \shuffle)$ generated by the letters
$w_i\ (1\le i \le r)$.

\subsubsection{}
Let $G$ be a simply connected complex Lie group with Lie algebra $\g$,
and let $N$ be a maximal unipotent subgroup of $G$ with 
Lie algebra $\nn$.
It is known that $U^*_\C$ is isomorphic to the algebra
$\C[N]$ of regular functions on $N$.
Hence the specialization $\varphi$ of $\Phi$ at $q=1$ may be regarded as an 
embedding of $\C[N]$ in the shuffle algebra $(\F_\C, \shuffle)$. 
The next proposition gives a direct description of $\varphi$.
Let $\underline{e}_i\ (1\le i\le r)$ be the Chevalley generators
of the Lie algebra $\nn$. Denote by
\[
x_i(t) = \exp(t\underline{e}_i), \qquad (1\le i \le r,\ t\in\C),
\]
the corresponding root subgroups in $N$.
\begin{proposition}\label{varphi}
Let $f\in\C[N]$ be homogeneous of degree $\nu\in Q^+$.
Let $w=w[i_1,\ldots ,i_k]\in\M$ be of degree $\nu$.
The coefficient of the word $w$ in the $\M$-expansion of $\varphi(f)$ 
is equal to the coefficient of the monomial $t_1\cdots t_k$ in the polynomial function
\[
(t_1,\ldots ,t_k) \mapsto f(x_{i_1}(t_1)\cdots x_{i_k}(t_k))\,.
\] 
Equivalently, we have
\[
\varphi(f) = \sum_{w=w[i_1,\ldots ,i_k],\ |w|=\nu} \left(
{\partial\over\partial t_1}\cdots{\partial\over\partial t_k}
f(x_{i_1}(t_1)\cdots x_{i_k}(t_k))|_{t_1=\cdots =t_k=0} \right)w\,.
\]
\end{proposition}
\proof
Let $u\in U^*_\A$ be homogeneous of degree $\nu$, and let
$f\in\C[N]$ be its specialization at $q=1$.
The group $N$ acts on $\C[N]$ by right translations:
\[
(x_i(t)f)(x) = f(xx_i(t)),\qquad (x\in N,\ t\in\C,\ 1\le i\le r)).
\]
Accordingly, the Lie algebra acts via the infinitesimal right translation
operators 
\[
\underline{e}_i(f)(x)={d\over dt} f(xx_i(t)) |_{t=0}.
\]
These are the specializations at $q=1$ of the endomorphisms
$e'_i\in\End U_q(\nn)$.
In particular, for any $i$, $\underline{e}_i(f)$ is homogeneous of degree $\nu-\a_i$.
It is easy to check that
\[
\underline{e}_{i_1}\cdots \underline{e}_{i_k}(f)(x) = 
{\partial\over\partial t_1}\cdots{\partial\over\partial t_k}
f(xx_{i_1}(t_1)\cdots x_{i_k}(t_k))|_{t_1=\cdots =t_k=0}
\,.
\]
If $\nu=\a_{i_1}+\cdots +\a_{i_k}$ this function of $x$ is a constant, equal
to the specialization at $q=1$ of $e'_{i_1}\cdots e'_{i_k}(u)$,
that is, to the coefficient of $w$ in $\varphi(f)$,
and the proposition follows.
\cqfd

One can also describe the inverse map 
$\varphi^{-1} : \U^*_\C \longrightarrow \C[N]$.
\begin{proposition}\label{phi_inverse}
Let $u=\sum_{w\in\M} \gamma(w) \, w$ be a homogeneous element
of $\U^*_\C$ of degree $\nu$, and let $f=\varphi^{-1}(u)$.
For any $(i_1,\ldots ,i_k)\in [1,r]^k$
and $(t_1,\ldots,t_k)\in\C^k$ we have
\[
f(x_{i_1}(t_1)\cdots x_{i_k}(t_k)) = 
\sum_w \gamma(w) 
{t_1^{a_1}\cdots t_k^{a_k}\over a_1!\cdots a_k!} 
\]
where the sum is over all $w =w_{i_1}^{a_1}\cdots w_{i_k}^{a_k}$ such
that $a_1\a_{i_1}+\cdots+ a_k\a_{i_k} = \nu$. 
\end{proposition}
\proof
This follows immediately from the identity
\[
f(xx_i(t))=(x_i(t)f)(x)=((\exp t \underline{e}_i)f)(x)
= \sum_{k\ge 0} {t^k\over k!} (\underline{e}_i^k f)(x)
\]
and the proof of the previous proposition.
\cqfd


\section{Good words and monomial bases} \label{SECT2}

From now on, we fix an arbitrary total order on the set
$\Pi=\{\a_1,\ldots ,\a_r\}$
of simple roots of~$\g$.
The alphabet $I=\{w_1,\ldots ,w_r\}$ is given the corresponding total 
order, and $\M$ the associated lexicographic order.
All these orders will be denoted by $<$.

\subsection{}
To $w=w[i_1,\ldots,i_k]$ we associate 
$D_w:=\ep_{i_1}\cdots \ep_{i_k} \in \End\,\F$.
(For $w=w[\,]$ we set $D_w = \Id_{\F}$.)
We have
\begin{equation}\label{partD}
D_w(\Phi(u))=\Phi(\D_w(u)),\qquad (u\in U_q(\nn)).
\end{equation}
For a homogeneous element $f$ of $\F$ we denote by $\max(f)$ the 
maximal word $w\in\M$
such that $|w|=|f|$ and $D_w(f)\not = 0$, that is, the largest word
occuring in the expansion of $f$.

\begin{definition}
A word $w\in\M$ is called good if there exists a homogeneous $u\in\U$ such that
$w=\max(u)$.
\end{definition}
The set of good words is denoted by $\G$.
Good words have been introduced by Lalonde and Ram for Lie algebras
and universal enveloping algebras \cite{LR}, and used by Rosso
in the context of quantum groups \cite{R3}.
Note that our definition is different from that of \cite{LR,R3}.
It will be shown in Lemma~\ref{charactgood} that the two definitions
are equivalent.
 
\begin{proposition}\label{monbasis} 
{\rm (i)}\ There is a unique basis of homogeneous vectors 
$\{m_g\ |\ g\in \G\}$ of $\U$ such that 
\[
D_{g_1}(m_{g_2})=\delta_{g_1g_2}, \qquad (g_1,g_2\in\G,\ |g_1|=|g_2|).
\]  
{\rm (ii)}\ $\{e_g\ |\ g\in \G\}$ is a basis of $U_q(\nn)$.
\end{proposition}
\proof
Let $\nu\in Q^+$. 
Let $\U_\nu$ be the homogeneous component of degree $\nu$ of $\U$,
$B_\nu$ a basis of $\U_\nu$, and $\{g_1,\ldots ,g_m\}$
the subset of $\G$ consisting of all words of weight $\nu$ arranged in
increasing order.
There exists at least one element of $B_\nu$, say $b_m$, such
that $\max(b_m)=g_m$.
By rescaling $b_m$ and subtracting appropriate multiples of it 
from the other elements of $B_\nu$ we can arrange that 
$g_m$ appears in $b_m$ with coefficient 1 and does not occur in any other
vector of $B_\nu$.
(Here we abuse notation and still denote by $B_\nu$ the basis
obtained after these operations.) 
Similarly, there exists a vector in $B_\nu \setminus \{b_m\}$,
say $b_{m-1}$, such that $\max(b_{m-1})=g_{m-1}$, and we can modify
$B_\nu$ in such a way that $D_{g_{m-1}}(b_{m-1})=1$ and $g_{m-1}$ 
occurs in no other element than $b_{m-1}$.
Repeating this process we get a subset $\{b_1,\ldots ,b_m\}$
of $B_\nu$ such that $D_{g_i}(b_j)=\delta_{ij}$.
Finally, $B_\nu = \{b_1,\ldots ,b_m\}$,
since otherwise there would be some $b\in B_\nu$ with 
$\max(b)\not = g_i$ for all $i$, which is impossible.
Proceeding in the same way in every weight space of $\U$ we
obtain a basis $\{m_g\ |\ g\in\G\}$ as in (i).
The unicity is clear. 

By (\ref{eq:1.2}) and (\ref{partD}) we see that the basis of
$U_q(\nn)$ adjoint to $\{\Phi^{-1}(m_g)\ |\ g\in\G\}$ is  
$\{e_{\tau(g)}\ |\ g\in\G\}$, where $\tau$ is as in Proposition~\ref{proptau}.
Finally, applying the anti-automorphism of $U_q(\nn)$ which fixes the 
generators $e_i$, we obtain that $\{e_g\ |\ g\in\G\}$ is a basis of
$U_q(\nn)$.
\cqfd

\begin{lemma}[\cite{LR}]\label{factor} 
Every factor of a good word is good.
\end{lemma}
\proof
Let $w=w[i_1,\ldots ,i_k]\in\G$ and let $u\in\U$ be such that $w=\max(u)$.
One checks easily that 
$
w[i_1,\ldots,i_j] = \max(\ep_{i_{j+1}}\cdots \ep_{i_k}(w)),
\ (1\le j<k).
$
We know that $\U$ is stable under the endomorphisms $\ep_i$, hence
all left factors of $w$ are good.
To conclude in the general case, we may introduce the endomorphisms
$\ep_i^\dag$ of $\F$ defined by
$
\ep_i^\dag(w[i_1,\ldots,i_k]) = \delta_{i,i_1} w[i_2,\ldots ,i_k].
$
In other words, $\ep_i^\dag=\tau \circ \ep_i \circ \tau$, where $\tau$
is as in Proposition~\ref{proptau}.
This shows that $\U$ is stable under $\ep_i^\dag$. 
Therefore,
$ 
w[i_h,\ldots,i_j] = \max(\ep_{i_{j+1}}\cdots
\ep_{i_k}\ep_{i_{h-1}}^\dag\cdots \ep_{i_1}^\dag(u))
$
is good for all $1<h\le j<k$.
\cqfd  

We are now going to study the set $\G$ of good words, and find an 
explicit description of it in terms of Lyndon words.

\subsection{}\label{SectLynd}
A word $l=w[i_1,\ldots,i_k]\in\M$ is a Lyndon word if it is smaller 
than all its proper right factors, that is,
\[
l < w[i_j,\ldots,i_k], \qquad (j=2,\ldots,k).  
\]
We shall denote by $\LL$ the set of Lyndon words in $\M$. 
For properties of Lyndon words which are not proved here,
see \cite{Lo} chapter 5, or \cite{Re}. 

We have the following inductive characterization of Lyndon words,
namely, $l$ is a Lyndon word if and only if $l\in I$ or
$l$ has a non trivial factorization $l=l_1l_2$ where 
$l_1$ and $l_2$ are Lyndon words and $l_1<l_2$.

For $l\in\LL\setminus I$, write $l=l_1l_2$ with 
$l_2$ a Lyndon word of maximal length.
It is known that $l_1$ is then also a Lyndon word,
and $l=l_1l_2$ is called the standard factorization of $l$.

Similarly, write $l=l_1^*l_2^*$ where $l_1^*$ is a Lyndon
word of maximal length. Then $l_2^*$ is also a Lyndon word
and we shall call $l=l_1^*l_2^*$ the co-standard factorization
of $l$.
This follows from the next lemma which gives a description
of $l_2^*$.

\begin{lemma}\label{costandard}
Let $l=l_1^*l_2^*$ be the co-standard factorization of
$l\in\LL$.
Then $l_2^*$ is of the form
\[
l_2^*=(l_1^*)^k f x\,,
\]
where $k\in\N$, $f$ is a left factor of $l_1^*$ (possibly empty),
and $x$ is a letter such that $fx>l_1^*$.
\end{lemma}
\proof
Let $m$ be a non trivial left factor of $l_2^*$. 
We want to prove that $m=(l_1^*)^k\,f$ for some $k\in\N$ 
and some left factor $f$ of $l_1^*$.
We will proceed by induction on $\ell(m)$.

Note first that by definition of the co-standard factorization,
$l_1^*m$ is not a Lyndon word. Thus $l_1^*m$ has a right factor
$\le l_1^*m$. This factor cannot be of the form $dm$ for some 
right factor $d$ of $l_1^*$, since $d>l_1^*$ and $\ell(d)<\ell(l_1^*)$
imply $dm > l_1^*m$.
Therefore this factor is a right factor of $m$.
In particular, if $\ell(m)=1$, we obtain that $m$ is less or
equal to the first letter $a$ of $l_1^*$, and since $l_1^*l_2^*$
is a Lyndon word, we must have $m=a$, which proves the claim in this
case.

Suppose now that $\ell(m)>1$ and write $m=m'y$ where $y$ is a letter.
By induction we may assume that $m'=(l_1^*)^k f$ for some $k\in\N$ 
and some non trivial left factor $f$ of $l_1^*$, and we have to prove that
$fy$ is a left factor of $l_1^*$ (possibly equal to $l_1^*$).
There exists a right factor $d$ of $m$ such that $d\le l_1^*m$.
We have $d=d'(l_1^*)^l\,fy$ for some $l\le k$ and some right factor
$d'$ of $l_1^*$. In fact $d'$ must be empty, otherwise,
since $l_1^*$ is Lyndon we would have $d'>l_1^*$ and 
$\ell(d')<\ell(l_1^*)$, hence 
\[
d=d'(l_1^*)^l\,fy > l_1^*(l_1^*)^k\,fy=l_1^*m\,.
\]
Therefore, $d=(l_1^*)^lfy\le (l_1^*)^{k+1}fy$
for some $l\le k$.
It follows that $fy\le (l_1^*)^i\,fy$ for some $i>0$,
and since $\ell(fy)\le \ell(l_1^*)$ we have in fact
$fy\le l_1^*$. Now either $fx$ is a left factor of $l_1^*$
or there is a $j\le \ell(fy)$ such that the $j$th letter
of $fy$ is strictly smaller than the $j$th letter of $l_1^*$.
The second case is impossible since then $l_1^*l_2^*$ 
could not be a Lyndon word. Therefore $fy$ is a left factor
of $l_1^*$.

Now we can apply this to the longest strict left factor $m$
of $l_2^*$, and we obtain that $l_2^*$ is as we claimed.
Moreover, since $l_1^*l_2^*$ is a Lyndon word, the letter $x$
has to be such that $fx>l_1^*$.
\cqfd

\begin{lemma}\label{lemmaLyndShuf}
Let $l \in \LL$ and $w\in\M$ with $l\ge w$. The largest word occuring
in the shuffle product of $l$ and $w$ is $lw$.
\end{lemma}
\proof
Let $u=u_1\cdots u_s$ be the largest word occuring
in the shuffle product of $l=a_1\cdots a_r$ and $w=b_1\cdots b_t$.
Suppose that $u\not = lw$, and suppose that the first
letter in which these two words differ is $u_k$.
Since $u$ occurs in the shuffle product, $k$ has to be
less or equal to the length $r$ of $l$, and $u_k=b_1>a_k$. 
Hence $a_k < w \le l$. Since $l$ is a Lyndon word, its 
smallest letter is its first letter, so this would force
$a_k=a_1$, hence $a_1<b_1$ and $l<w$, contrary to our
assumption.
\cqfd

Let $\GL$ denote the subset of $\G$ consisting of all good Lyndon words.
\begin{proposition}\label{GLG}
Let $l\in\GL$ and $g\in\G$ with $l\ge g$.
Then $lg\in\G$.
\end{proposition}
\proof
Let $u,v$ be homogeneous elements of $\U$ such that
$\max(u)=l$ and $\max(v)=g$.
Rescaling $u$ and $v$ if necessary we can assume that 
$u=l+r$ and $v=g+s$ where $r$ (\resp $s$) is a linear
combination of words $<l$ (\resp $<s$).
We have
$
u\shuf v = l\shuf g + l\shuf s + r\shuf g + r\shuf s.
$
By Lemma~\ref{lemmaLyndShuf}, $\max(l\shuf g) = lg$.
Now if $w$ and $w'$ are words such that $|w|=|l|$,
$|w'|=|g|$, $w\le l$ and $w'\le g$, any word occuring
in the shuffle of $w$ and $w'$ will be less or equal
to the corresponding word in the shuffle of $l$ and $g$,
so $\max(u\shuf v) = lg$.
\cqfd 

It is well known \cite{Lo,Re} that
every word $w\in\M$ has a unique factorization
$w = l_1\cdots l_k$
where $l_1,\ldots ,l_k \in \LL$ and $l_1\ge \cdots \ge l_k$.

\begin{proposition}[\cite{LR}]\label{factgood}
A word $g$ is good if and only if it is of the form
\[
g = l_1\cdots l_k,\qquad l_1\ge \cdots \ge l_k,
\]
where $l_1, \ldots , l_k$ are good Lyndon words.
\end{proposition}
\proof
By Lemma~\ref{factor}, if $g$ is good, its canonical factorization
as a non-increasing product of Lyndon words has good factors.
The converse follows immediately from Proposition~\ref{GLG}.
\cqfd

\begin{proposition}[\cite{LR}]
The map $l\mapsto |l|$ is a bijection from $\GL$ to $\Delta^+$.
\end{proposition}
\proof
By Proposition~\ref{monbasis} and Proposition~\ref{factgood}
the products
\[
e_{l_1}^{n_1}\cdots e_{l_k}^{n_k},\qquad 
(n_1,\ldots n_k \in\N,\ l_1,\ldots,l_k\in\GL,\ l_1>\cdots >l_k) 
\]
form a basis of $U_q(\nn)$.
This implies that the generating series
of the dimensions of the homogeneous components of $U_q(\nn)$ is
equal to
\[
\sum_{\nu \in Q^+} \dim U_q(\nn)_\nu \, \exp\nu = \prod_{l\in \GL}
{1\over 1-\exp|l|}\,.
\]
On the other hand it is well-known that
\[
\sum_{\nu \in Q^+} \dim U_q(\nn)_\nu \, \exp\nu = \prod_{\b\in\Delta^+}
{1\over 1-\exp\b}\,,
\]
and by comparing the two expressions the claim follows.
\cqfd

We shall denote by $\b\mapsto l(\b)$ the inverse of the above 
bijection. It is an embedding of $\Delta^+$ in~$\LL$.
We will call it a Lyndon covering of $\Delta^+$.


\section{PBW-type bases}\label{SECT2bis}
In this section, we introduce following Lalonde-Ram \cite{LR} and 
Rosso \cite{R3} another basis $\{r_g\}$ of $\U$ labelled by good
words, the Lyndon basis.
Then we show that this basis is up to normalization the image 
under the anti-automorphism $\si$ of
a basis $\{E_g\}$ of PBW-type, as defined by Lusztig \cite{L}.
This allows us to prove that $\max(E_g)=g$ (Theorem~\ref{theomax}).
In \ref{LyndCover} we also provide an algorithm for computing explicitly
the map $\b\mapsto l(\b)$.
This works for any root system and any ordering of the simple
roots, and is simpler than the procedure of \cite{LR} which needs
some case-by-case discussion. 
Finally, we prove that the normalization coefficient $\k_g$ between
$\si(r_g)$ and $E_g$ is a bar-invariant Laurent polynomial
(Proposition~\ref{propkappa}), which will be used in \ref{newproof}.

\subsection{} \label{subsec:3.2}
For homogeneous elements $f_1,f_2\in\F$ we define
\[
[f_1,f_2]_q:=f_1f_2-q^{(|f_1|,|f_2|)}\,f_2f_1.
\]
Let $l\in\LL$.
We define inductively the $q$-bracketing $[l]\in\F$
by $[l]=l$ if $l$ is a letter, and otherwise
$[l] = [[l_1],[l_2]]_q$,
where $l=l_1l_2$ is the co-standard factorization of $l$.
\begin{proposition}\label{PropRosso0}
$[l]= l + r$ where $r$ is a linear
combination of words $>l$. 
\end{proposition}
\proof
We argue by induction on the length of $l$.
If $l$ is a letter, the statement is obvious. Otherwise 
$[l] = [[l_1],[l_2]]_q$ and we can assume by induction
that $[l_1]= l_1 + r_1$ and $[l_2]= l_2 + r_2$ where $r_1$ and $r_2$
are linear combinations of words $>l_1$ and $>l_2$ respectively.
Hence, 
\[
[[l_1],[l_2]]_q=[l_1,l_2]_q+[r_1,l_2]_q+[l_1,r_2]_q+[r_1,r_2]_q. 
\]
The first bracket is $l_1l_2-q^{(|l_1|,|l_2|)}\,l_2l_1$, and 
since $l_1l_2$ is a Lyndon word, $l_2l_1>l_1l_2$.
Clearly, all words occuring in the other brackets are 
either $>l_1l_2$ or $>l_2l_1>l_1l_2$, and the statement follows.
\cqfd

Let $w=l_1\cdots l_k$ be the canonical factorisation of $w$
as a non-increasing product of Lyndon words.
We define 
$[w]:=[l_1]\cdots [l_k]\in\F.$
\begin{proposition}\label{prop8}
$\{[w]\ |\ w\in\M\}$ is a basis of $\F$.
\end{proposition}
\proof
It follows easily from Proposition~\ref{PropRosso0} 
that $[w]=w+s$ where $s$ is a linear combination of words $>w$.
Hence the transition matrix from the basis $\{w\}$ to the
family of vectors $\{[w]\}$ is unitriangular.
\cqfd

\subsection{}
Let $\Xi$ be the algebra homomorphism from $(\F,\cdot)$ to
$(\F,\shuf)$ such that $\Xi(w_i) = w_i$ for every letter $w_i$,
that is, each word is mapped by $\Xi$ to the quantum shuffle
product of its letters.
Clearly, $\Xi(\F)=\U$.

\begin{lemma}\label{charactgood}
The word $w$ is good if and only if it cannot be expressed
modulo $\ker \Xi$ as a linear combination of words $v>w$. 
\end{lemma}
\proof
Write $w=w[i_1,\ldots,i_k]$.
Suppose that $w$ can be expressed modulo $\ker \Xi$ as a linear 
combination of words $v>w$, that is, there exists a relation of the form
\begin{equation}\label{relation}
w_{i_1}\shuf \cdots \shuf w_{i_k} =
\sum_{v=w[j_1,\ldots,j_k]>w} x_v\, w_{j_1}\shuf \cdots \shuf w_{j_k}
\end{equation}
for some scalars $x_v\in\Q(q)$.
Using the isomorphism $\Phi^{-1} : \U \longrightarrow U_q(\nn)$
this is equivalent to
\[
e_{i_1} \cdots e_{i_k} =
\sum_{v=w[j_1,\ldots,j_k]>w} x_v\, e_{j_1} \cdots e_{j_k},
\]
and since the algebra generated by the $e'_i$ is isomorphic to $U_q(\nn)$ 
this is in turn equivalent to
\[
\D_w =
\sum_{v>w} x_v\,\D_v.
\]
Therefore, if for some homogeneous $u\in \U$ of weight $|u|=|w|$
one has $D_w(u)\not = 0$, then there exists a $v>w$ such that
$D_v(u)\not = 0$, and $w\not =\max(u)$.
Hence $w$ is not good.
Let us denote by $\HHH$ the set of words $w$ which satisfy
no relation of the form (\ref{relation}).
We have proved that $\G \subset \HHH$.

Conversely, it is easy to prove that $\{e_w \ | \ w\in\HHH\}$
is a basis of $U_q(\nn)$. 
Indeed, this set contains the monomial basis of
Proposition~\ref{monbasis},
and it is linearly independent, since if we had a linear 
relation between words of $\HHH$ we could express the smallest
one in terms of the others and it would not belong to $\HHH$.
Hence $\G=\HHH$, as required.
\cqfd

Note that in \cite{LR,R3}, Lemma~\ref{charactgood} is taken as 
the definition of a good word.

For $g\in\G$, let us write $r_g=\Xi([g])$.
\begin{proposition}[\cite{R3}]\label{PropRosso}
$\{r_g\ |\ g\in\G\}$ is a basis of $\U$.
\end{proposition}
\proof
Note that for any word $w$, we have $\Xi(w)=\Phi(e_w)$.
As in the proof of Proposition~\ref{prop8}, for 
$g\in\G$ we have $[g] = g + \sum_{w>g}x_{gw}\, w$.
Thus $r_g = \Phi(e_g) + \sum_{w>g} x_{gw}\, \Phi(e_w)$. 
By Lemma~\ref{charactgood}, this last sum can be rewritten as
$
r_g = \Phi(e_g) + \sum_{h>g} y_{gh}\, \Phi(e_h),
$
where the words $h$ are good. 
Hence, the transition matrix from the basis $\Phi(\{e_g\ |\ g\in\G\})$
to $\{r_g\ |\ g\in\G\}$ is unitriangular.  
\cqfd

We call $\{r_g\ |\ g\in\G\}$ the Lyndon basis of $\U$.
\begin{theorem}[\cite{R3}]\label{thRosso1}
The Lyndon basis has the following form
\[
\{r_{l_1}\shuf\cdots\shuf r_{l_k}\ | \ k\in\N,\ l_1,\ldots ,
 l_k\in\GL,\ l_1\ge \cdots \ge l_k\}\,.
\]
\end{theorem}
\proof
By definition of $[g]$, if $g=l_1\cdots l_k$ is the canonical factorization
of $g$ as a non-increasing product of Lyndon words, 
we have $r_g=r_{l_1}\shuf\cdots\shuf r_{l_k}$, and by
Lemma~\ref{factor} each factor $l_k$ is good.
Conversely, if $l_1,\ldots , l_k$ are good Lyndon words 
and $l_1\ge \cdots \ge l_k$ then by Proposition~\ref{factgood},
$g=l_1\cdots l_k$ is good. 
\cqfd

\begin{proposition}\label{major}
Let $\b_1,\b_2\in\Delta^+$ be such that $\b_1+\b_2=\b\in\Delta^+$
and $l(\b_1)<l(\b_2)$.
Then $l(\b_1)l(\b_2)\le l(\b)$.
\end{proposition}
\proof
As seen in the proof of Proposition~\ref{PropRosso}, the transition
matrix from the basis $\{\Xi(g)\ |\ g\in\G\}$ to the basis
$\{r_g\ |\ g\in\G\}$ is unitriangular. 
Hence, writing $l_1=l(\b_1)$ and $l_2=l(\b_2)$, we have 
\begin{eqnarray*}
r_{l_1}\shuf r_{l_2}& =&
\left(\Xi(l_1)+\sum_{h_1>l_1;\,h_1\in\G} y_{l_1h_1}\,\Xi(h_1)\right)
\shuf
\left(\Xi(l_2)+\sum_{h_2>l_2;\,h_2\in\G} y_{l_2h_2}\,\Xi(h_2)\right)\\
&=&\sum_{g\ge l_1l_2} z_g \,r_g\,
\end{eqnarray*}
for some $z_g\in\Z[q,q^{-1}]$.
Recall from \ref{special2} the $\A$-subalgebra $U_\A$ of $U_q(\nn)$.
Let $x\mapsto \underline{x}$ denote the specialization $q\mapsto 1$
from $U_\A$ to $U(\nn)$.
For $l\in\GL$ set $s_l=\Phi^{-1}(r_l)$.
Then $\underline{s_l}\in\nn$
(this is an iterated bracket of Chevalley generators 
$\underline{e_i}$), 
and $[\underline{s_{l_1}},\underline{s_{l_2}}]$
belongs to the weight space of weight $\b$ of $\nn$. 
By hypothesis, this weight space is $1$-dimensional and spanned
by $\underline{s_{l(\b)}}$.
Hence,
\[ 
\underline{s_{l_1}}\,\,\underline{s_{l_2}}=\underline{s_{l_2}}\,\,\underline{s_{l_1}}
+c\,\underline{s_{l(\b)}}
\]
for some $c\in\Z^*$.
It follows that 
$z_{l(\b)}\not = 0$, which implies that $l(\b)\ge l_1l_2$.
\cqfd

\subsection{}\label{LyndCover}
Proposition~\ref{major} implies the following simple inductive rule
for determining the set $\GL$ of good Lyndon words.
If $\b=\a_i$ is a simple root, then $l(\b)=w_i$.
If $\b$ is not a simple root 
there exists a factorization $l(\b)=l_1l_2$ with $l_1$ and $l_2$ Lyndon.
By Lemma~\ref{factor}, $l_1=l(\b_1)$ and $l_2=l(\b_2)$ for some
$\b_1,\,\b_2\in\Delta^+$.
By induction
we may assume that we know $l(\gamma)$
for any $\gamma\in\Delta^+$ of height smaller than the height of $\b$.
Let 
\[
C(\b)=\{(\b_1,\b_2)\in\Delta^+\times\Delta^+\ |\ \b_1+\b_2=\b,\ l(\b_1)<l(\b_2)\}.
\]
Then, by Proposition~\ref{major}, we get
\begin{proposition} 
$\quad l(\b)=\max\{l(\b_1)l(\b_2)\ |\ (\b_1,\b_2)\in C(\b)\}.$
\cqfd
\end{proposition}

Note that the sets $\LL$, $\G$, $\GL$ depend on the choice of a 
total order on $\Pi$, and that we have $r!$ possible choices.
In \cite{LR} the sets $\GL$ are calculated for all root systems
and for a particular total order on $\Pi$
(see also section~\ref{SECT6}).

By Proposition~\ref{factgood}, we can calculate the set $\G$ of good
words by taking the non-increasing products of elements of $\GL$. 
Note that, by Proposition~\ref{monbasis},
we have thus obtained for each total order on $\Pi$ a simple
and explicit monomial basis of $U_q(\nn)$.
This basis seems to be different from the monomial bases of 
Chari and Xi \cite{CX} and Reineke \cite{Rei}. 

\subsection{}
Since $\GL$ is totally ordered (lexicographically) we obtain
a total order (still denoted by $<$) on $\Delta^+$.
In \cite{R3} the following key fact is stated.
\begin{proposition}[\cite{R3}]\label{keylemma}
The order $<$ on $\Delta^+$ is convex, that is, 
if $\b_1$ and $\b_2$ are elements of $\Delta^+$
such that $\b_1+\b_2=\b$ belongs to $\Delta^+$, 
then $\b_1<\b<\b_2$ or $\b_2<\b<\b_1$.
\end{proposition}
Note that by Proposition~\ref{major}, we have that if
$l(\b_1)<l(\b_2)$ and $\b_1+\b_2=\b$ then 
\[
l(\b)\ge l(\b_1)l(\b_2) > l(\b_1).
\]
On the other hand, if $l(\b)=l(\b_1)l(\b_2)$, since
$l(\b)$ is a Lyndon word, $l(\b)<l(\b_2)$.
It only remains to prove that, even when $l(\b)>l(\b_1)l(\b_2)$
we have $l(\b)<l(\b_2)$.
\begin{corollary}\label{smallest}
Let $\b\in\Delta^+$.
The good Lyndon word $l$ of weight $\b$ is the 
smallest good word of weight $\b$.
\end{corollary}
\proof
Let $g\not = l$ be a good word of weight $\b$ and let
$g=l_1\cdots l_k$ be its unique expression
as a non-increasing product of good Lyndon words.
Let $\b_i=|l_i|\ (1\le i \le k)$.
If $g<l$ then $l_1<l$. 
Indeed, Melan\c con has shown that if 
$w=m_1\cdots m_r$ and $w'=m'_1\cdots m'_s$
are the factorizations into non-increasing products
of good Lyndon words of $w$ and $w'$,
we have $w>w'$ if and only if there exist $j$
such that $m_i=m'_i$ for $i<j$ and $m_j>m'_j$ \cite{Me}.
Therefore $l_i < l$ for all $i=1,\ldots ,k$,
hence we have $\b_1+\cdots +\b_k = \b$ with
all $\b_i<\b$, contrary to the fact that $<$
is convex.
\cqfd 

\subsection{}\label{PBWmain}
It is well-known \cite{P} that any convex ordering 
$\b_1<\cdots <\b_n$ of $\Delta^+$
arises from a unique reduced decomposition $w_0=s_{i_1}\cdots s_{i_n}$ 
of the longest element of the Weyl group $W$ in the following
way:
\[
\b_1=\a_{i_1},\quad \b_2=s_{i_1}(\a_{i_2}) \quad ,\ldots ,\quad  
\b_n=s_{i_1}\cdots s_{i_{n-1}}(\a_{i_n})\,.
\]
To this data Lusztig associates a PBW-type basis of $U_q(\nn)$ 
\[
E^{(a_1)}(\b_1)\cdots E^{(a_n)}(\b_n),\qquad (a_1,\ldots ,a_n)\in
\N^n,
\] 
defined using the braid group action on $U_q(\nn)$ (\cite{L} 40.2.2).
(We choose the action via the operators $T'_{i,-1}$ of
\cite{L} 37.1.3, with $q=v^{-1}$.) 
Let us fix from now on the PBW-type basis associated with
the convex ordering on $\Delta^+$ coming from its Lyndon covering
$\GL$.
\begin{theorem}\label{RossoPBW}
For all $\b\in\Delta^+$, the vectors $\Phi(E(\b))$ and 
$r_{l(\b)}=\Xi([l(\b)])$
are proportional.
\end{theorem}
\proof
We argue by induction on the height $k$ of $\b$.
If $\b$ is a simple root, the claim is trivial.
Suppose that $k>1$ and that the result is proved for all 
roots of height $\le k-1$. 
We can write $\b$ as a sum $\b_1+\b_2$ of two positive
roots $\b_1 < \b_2$ (\cite{Bo}, Prop. 19), and clearly these roots 
both have height $\le k-1$.
Among all such decompositions, pick up the one for which $\b_1$
is maximum and denote it by $\b = \b^*_1+\b^*_2$.
By a result of Levendorskii and Soibelman \cite{LS} 
(see also \cite{CP}, 9.3),
\[ 
E(\b^*_1)E(\b^*_2)-q^{(\b^*_1,\b^*_2)}E(\b^*_2)E(\b^*_1)
\]
is a linear combination of products 
$
E(\b_{i_1})\cdots E(\b_{i_s})
$
where
$\b_{i_1}+\cdots +\b_{i_s} = \b$
and $\b^*_1<\b_{i_j}<\b^*_2$ for every $j$.
(Note that the $\b_{i_j}$ are not necessarily distinct.)
Suppose there occurs in this linear combination a term
other than $E(\b)$, that is, a term with $s>1$, and let us consider it.
By \cite{Bo}, Prop. 19, we have (after renumbering the $\b_{i_j}$'s 
if necessary) that $\b_{i_1}+\cdots +\b_{i_j}\in\Delta^+$ for
every $j$.
In particular $\b'=\b_{i_1}+\cdots +\b_{i_{s-1}}$ and $\b''=\b_{i_s}$
are two positive roots with $\b'+\b''=\b$.
Therefore, by definition of $\b_1^*$
and because $<$ is convex,
either $\b'<\b_1^*<\b<\b''$ or $\b''<\b_1^*<\b<\b'$.
In the second case we would have $\b_{i_s}<\b^*_1$ which is impossible.
The first case is also impossible since if all $\b_{i_j}>\b^*_1$
then $\b'>\b^*_1$ by convexity.
Hence, 
\[
E(\b^*_1)E(\b^*_2)-q^{(\b^*_1,\b^*_2)}E(\b^*_2)E(\b^*_1)
\]
is proportional to $E(\b)$.

On the other hand, let us consider the element
\[
r=r_{l(\b^*_1)}r_{l(\b^*_2)}-
q^{(\b^*_1,\b^*_2)}r_{l(\b^*_2)}r_{l(\b^*_1)}.
\]
Let $l=l_1^*l_2^*$ be the co-standard factorization of $l=l(\b)$.
We have $l(\b_1^*)\ge l_1^*$ by definition of $\b_1^*$, and
$l(\b_1^*)l(\b_2^*)\le l_1^*l_2^*$ by Proposition~\ref{major}. 

This implies that $l(\b_1^*)=l_1^*$ and $l(\b_2^*)=l_2^*$.
Indeed, the two inequalities imply that $l_1^*$ is a left
factor of $l(\b_1^*)$, that is, $l(\b_1^*)=l_1^*m$.
Suppose that $m$ is not the empty word.
If $l(\b_1^*)l(\b_2^*)= l_1^*l_2^*$, then $l_1^*l_2^*$ would
not be the co-standard factorization, so we must have
$l(\b_1^*)l(\b_2^*)< l_1^*l_2^*$.
By Lemma~\ref{costandard}, $l_2^*=(l_1^*)^k\,f\,x$ for
some $k\in\N$, some left factor $f$ of $l_1^*$ and some 
letter $x$ such that $fx>l_1^*$.
Since $l(\b_1^*)$ is a Lyndon word, $m>l_1^*$.
On the other hand we must have $l_1^*m=l(\b_1^*)<l_1^*l_2^*$,
hence $m<l_2^*$, and this implies
that $l_1^*$ is a left factor of $m$.
Since for any $k'$ and any left factor $f'$ of
$l_1^*$, $(l_1^*)^{k'}f'$ is not a Lyndon word, we see that
the only possibility is $m=l_2^*$, but then $l(\b^*_2)$
would be empty, a contradiction. Therefore $m$ is empty,
and $l(\b_1^*)=l_1^*$. Finally, since there is only one
good Lyndon word of weight $\b_2^*$, we also have $l_2^*=l(\b_2^*)$. 

Hence $r=r_{l(\b)}$ and
since by induction $r_{l(\b_1^*)}$ and $r_{l(\b_2^*)}$
are proportional to $\Phi(E(\b_1^*))$ and $\Phi(E(\b_2^*))$ respectively, 
we conclude that $r_{l(\b)}$ is proportional to $\Phi(E(\b))$.
\cqfd

In \cite{Ri}, Ringel has proved a result similar to
the above theorem for the PBW-bases of Lusztig associated to 
reduced words for $w_0$ adapted to an orientation of the Dynkin 
diagram of $\g$, that is, for those bases coming from the theory 
of Hall algebras. 
Note that the convex orderings of $\De^+$ coming from Lyndon
coverings are in general different from those coming from
Hall algebras, as shown by the next example. 
Hence Theorem~\ref{RossoPBW} is different from Ringel's 
result, since the PBW-bases involved are not the same. 

\begin{example}{\rm
Let $\g$ be of type $D_4$, with Dynkin diagram numbered as in 
\ref{TypeD}. 
The Lyndon covering associated to the order $\a_1<\a_2<\a_3<\a_4$
corresponds to the reduced decomposition
\[
w_0=s_1s_3s_2s_4s_3s_1s_4s_3s_2s_4s_3s_4
\]
which is adapted to no orientation of the Dynkin diagram of $\g$.
}
\end{example}

\subsection{}\label{SUBS46}
For a good word $g=l(\b_1)^{a_1}\cdots l(\b_k)^{a_k}$,
where $\b_1>\cdots >\b_k$ and $a_1,\ldots, a_k\in\N^*$,
we will denote by
\begin{equation}\label{PBW}
E_g:=\Phi(E(\b_k)^{(a_k)}\cdots E(\b_1)^{(a_1)}),
\end{equation}
the corresponding vector of the PBW-type basis of $\U$,
and we will write
\begin{equation}
r_l = \l_l\, E_l \qquad (l\in\GL).
\end{equation}
Note that in (\ref{PBW}), the factors are taken in the order opposite to
the order used for defining $r_g$.
Recall from \ref{automorph} the anti-automorphism $\si$.

\begin{proposition}\label{proport}
{\rm (i)}\ Let $l=w\in\GL$.
We have $\si(r_l)=(-1)^{\ell(l)-1}q^{-N(|l|)}\,r_l$.

\noindent
{\rm (ii)}\ 
The vectors $E_g$ and $\si(r_g)$ are proportional for all $g\in\G$.
\end{proposition}
\proof
Since $\si$ is an anti-automorphism and $E_l$ is proportional 
to $r_l$ for $l\in\GL$, we see that (ii) follows immediately
from (i).
Let $l=l_1l_2$ be the co-standard factorisation of $l$, so that
\[
r_l = r_{l_1}\shuf r_{l_2} - q^{(|l_1|,|l_2|)}\,r_{l_2}\shuf r_{l_1}.
\]
Let us assume that (i) holds for $l_1$ and $l_2$. Then
\[
\si(r_l) = \si(r_{l_2})\shuf \si(r_{l_1}) - 
q^{-(|l_1|,|l_2|)}\,\si(r_{l_1})\shuf \si(r_{l_2})
=
(-1)^{\ell(l_1)+\ell(l_2)-1}q^{-(|l_1|,|l_2|)-N(|l_1|)-N(|l_2|)}\,
r_l
,
\]
and the result follows from the equality
\[
-(|l_1|,|l_2|)-N(|l_1|)-N(|l_2|) 
= -\sum_{1\le s<t \le k} (\a_{i_s}\,,\a_{i_t})
= -N(|l|)\,.
\]
\cqfd

In the sequel, we will write 
\begin{equation}
\si(r_g) = \k_g\, E_g \qquad (g\in\G).
\end{equation}
Write $g=l_1^{a_1}\ldots l_k^{a_k}$ where $l_1>\cdots >l_k\in\GL$
and $a_1,\ldots ,a_k\in\N^*$.
Then $r_g=r_{l_1}^{a_1}\cdots r_{l_k}^{a_k}$,
while $E_g=E_{l_k}^{(a_k)}\cdots E_{l_1}^{(a_1)}$,
where for $l=l(\b)\in\GL$ we set 
$E_l^{(a)}=E_l^a/[a]_i!$ if $(\b,\b)=(\a_i,\a_i)$.
It follows that, writing $[a]_l:=[a]_i$,
\begin{equation}\label{kprod}
\k_g =\prod_{j=1}^k \k_{l_j}^{a_j}\,[a_j]_{l_j}!\, .
\end{equation}

\begin{proposition}
For $l\in\GL$ we have
$\si(E_l)=(-1)^{\ell(l)-1}q^{N(|l|)}E_l$.
\end{proposition}
\proof
Let $\b_1<\b_2<\cdots <\b_n$ be the convex ordering of $\Delta^+$
associated with $\GL$, and let $w_0=s_{i_1}\cdots s_{i_n}$ be the
corresponding reduced decomposition of $w_0$.
Suppose that $|l|=\b_k$.
Then we have 
$E_l=T_{i_1}\cdots T_{i_{k-1}}(e_{i_k})$.
Using \cite{L} 37.2.4, we see that 
$\si(E_l)=(-1)^{A(l)}q^{B(l)}\,E_l$
where
\[
A(l)={1\over d_{i_1}}(\a_{i_1}\,,s_{i_2}\cdots s_{i_{k-1}}(\a_{i_k}))
     +{1\over d_{i_2}}(\a_{i_2}\,,s_{i_3}\cdots s_{i_{k-1}}(\a_{i_k}))
     +\cdots
     +{1\over d_{i_{k-1}}}(\a_{i_{k-1}}\,,\a_{i_k})\,,
\]
and 
\[
B(l)=(\a_{i_1}\,,s_{i_2}\cdots s_{i_{k-1}}(\a_{i_k}))
     +(\a_{i_2}\,,s_{i_3}\cdots s_{i_{k-1}}(\a_{i_k}))
     +\cdots
     +(\a_{i_{k-1}}\,,\a_{i_k})\,.
\]
On the other hand, an elementary calculation gives
\begin{eqnarray*}
\b_k& =& \a_{i_k}-{1\over
  d_{i_{k-1}}}(\a_{i_{k-1}}\,,\a_{i_k})\,\a_{i_{k-1}}
 -{1\over d_{i_{k-2}}}(\a_{i_{k-2}}\,,s_{i_{k-1}}(\a_{i_k}))\,\a_{i_{k-2}}\\
&& -\cdots
 -{1\over d_{i_1}}(\a_{i_1}\,,s_{i_2}\cdots s_{i_{k-1}}(\a_{i_k}))\,\a_{i_1}
\,.
\end{eqnarray*}
Therefore, writing $\b_k = \sum_{i=1}^r c_i\a_i$, we see that
$A(l)=1-\sum_{i=1}^r c_i=1-\ell(l)$.
Finally,
\[
N(|l|)=N(\b_k)= {1\over 2}\left((\b_k\,,\b_k)-\sum_{i=1}^r c_i (\a_i,\a_i)\right)
       = {1\over 2}\left((\a_{i_k}\,,\a_{i_k})-\sum_{i=1}^r 2\, c_i d_i\right) 
       = B(l)\,.
\]

\begin{proposition}\label{propkappa}
{\rm (i)}\ For $l\in\GL$ we have 
$\k_l=\overline{\k_l}=(-1)^{\ell(l)-1}q^{-N(|l|)}\l_l$.

\noindent
{\rm (ii)} For $g\in\G$ we have $\k_g=\overline{\k_g}\in\Z[q,q^{-1}]$.
\end{proposition}
\proof
Since $r_l=\l_lE_l$ we have 
\[
\si(r_l)=(-1)^{\ell(l)-1}q^{-N(|l|)}r_l=(-1)^{\ell(l)-1}q^{-N(|l|)}\l_lE_l.
\]
On the other hand 
\[
\si(r_l)=\overline{\l_l}\si(E_l)=(-1)^{\ell(l)-1}q^{N(|l|)}\overline{\l_l}E_l.
\]
Hence 
\[
\k_l=(-1)^{\ell(l)-1}q^{-N(|l|)}\l_l=(-1)^{\ell(l)-1}q^{N(|l|)}\overline{\l_l}
=\overline{\k_l},
\]
which proves (i).

Let us prove (ii).
As before, write $\A=\Z[q,q^{-1}]$ and
consider the $\A$-subalgebra $\U_\A$ of $\U$ generated by 
the elements $\Phi(e_i^{(k)})$.
It is known that $\{E_g\}$ is an 
$\A$-basis of $\U_\A$ (\cite{L}, 41.1.4).
By construction, $r_l$ is an iterated $q$-commutator
of generators $e_i$, thus it
belongs to $\U_\A$, and since $\si$ preserves $\U_\A$,
we have $\k_l\in\A$.
Finally, it follows from (i) and Equation (\ref{kprod}) that 
$\overline{\k_g}=\k_g\in\A$.

\begin{lemma}\label{PBW2mon}
For $g\in\G$ we have
$
E_g=\sum_{h\in\G,\, h\ge g} \a_{gh}(q) \, \Phi(e_{\tau(h)}),
$
where 
$\a_{gg}=\k_g^{-1}$.
\end{lemma}
\proof
We have 
$
r_g = \sum_{h\ge g,\ h\in\G} y_{gh}\, \Phi(e_h)\,,
$
with $y_{gg}=1$.
The result then follows from the relation $\si(r_g) = \k_g \, E_g$.
\cqfd

\subsection{}
We endow $\U$ with the nondegenerate symmetric bilinear
form $(\cdot\,,\,\cdot)$ obtained by transporting 
Kashiwara's form on $U_q(\nn)$ to $\U$ via $\Phi$.
It is known that the PBW-type bases of Lusztig are orthogonal
(\cite{L}, 38.2.3).
More precisely we have for $g=l(\b_1)^{a_1}\ldots l(\b_n)^{a_n}$ 
and $h=l(\b_1)^{b_1}\ldots l(\b_n)^{b_n}$ where $a_1,\ldots
,a_n,b_1,\ldots ,b_n\in\N$, 
\begin{equation}\label{SP1}
(E_g\,,\,E_h)=\delta_{gh}
\prod_{j=1}^n {(E(\b_j)\,,\,E(\b_j))^{a_j}\over \{a_j\}_{(\b_j,\b_j)}!}
\,,
\end{equation}
where for $\b=\sum_{i=1}^r c_i \a_i\in\Delta^+$,
\begin{equation}\label{SP2}
(E(\b)\,,\,E(\b)) = {\prod_{i=1}^r(1-q^{(\a_i,\a_i)})^{c_i}\over 1-q^{(\b ,\b)}}
\end{equation}
and for $m,p\in\N$,
\begin{equation}
\{m\}_p! = \prod_{j=1}^m {1-q^{jp}\over 1-q^p}\,.
\end{equation}
Following Lusztig (\cite{L}, 1.2.10) let us define another symmetric
bilinear form $\{\cdot\,,\,\cdot\}$ by setting
\begin{equation}
\{u\,,v\}=\overline{(\overline{u}\,,\overline{v})},\qquad (u,v\in\U).
\end{equation}
\begin{lemma}
For homogeneous elements $u, v\in\U$ of weight $\nu$ we have
$
\{u\,,v\}=q^{N(\nu)}\,(u\,,\tau(v))\,.
$
\end{lemma}
\proof
It is enough to check the lemma when $u, v$ run through 
two bases of $\U$. 
Let us take 
$u=m_g,\ v=\Phi(e_h), \ (g,h\in\G).$
We have $\overline{v} = v$ and 
writing 
$
m_g = g + \sum_{w<g,\,w\not\in\G} \gamma_w(q)\,w,
$
by Proposition~\ref{proptau},
\[ 
\overline{u}= q^{-N(|g|)}\left(\tau(g) + \sum_{w<g,\,w\not\in\G}
\gamma_w(q^{-1})\,\tau(w)\right).
\]
Hence,
$\{u\,,v\}=q^{N(|g|)}\,\delta_{gh}=q^{N(|g|)}(u,\tau(v))$.
\cqfd

\begin{proposition}\label{orthRosso}
The basis $\{r_g\}$ is orthogonal with respect to $\{\cdot\,,\cdot\}$.
\end{proposition}
\proof
It is known (\cite{L} 1.2.8) that 
$(\tau(u)\,,\tau(v))=(u,v)$ for all $u, v \in \U$. 
Hence we have 
\[
\{r_g\,,\,r_h\}=\overline{(\tau(\overline{r_g})\,,\tau(\overline{r_h}))}
=\overline{(\si(r_g)\,,\si(r_h))}
=\overline{\k_g}\,\overline{\k_h}\, \overline{(E_g\,,E_h)}\,,
\quad
(g,h\in\G),
\]
and the result follows from the orthogonality of $\{E_g\}$
with respect to $(\cdot\,,\,\cdot)$.
\cqfd

\begin{theorem}\label{theomax}
For $g\in\G$, we have $\max(r_g)=\max(E_g)=g$.
\end{theorem}
\proof
By the proof of Proposition~\ref{PropRosso} we have
$
r_g = \sum_{k\ge g,\,k\in\G} y_{gk}\,\Phi(e_k)\,,\ (g\in\G),
$
where $y_{gg}=1$.
On the other hand 
$
\{\Phi(e_k)\,,m_h\} = q^{N(|k|)}\,(\Phi(e_{\tau(k)})\,,m_h)
=q^{N(|k|)}\,\delta_{kh}\,.
$
Hence, by Proposition~\ref{orthRosso}, we have
\[
m_h 
= \sum_{g\in\G} \{r_g\,,m_h\}\, {r_g\over \{r_g\,,r_g\}}
= \sum_{g\le h,\,g\in\G} q^{N(|g|)}\,y_{gh}\, {r_g\over \{r_g\,,r_g\}}\,,
\]
therefore
\[
r_h = q^{-N(|h|)}\,\{r_h\,,r_h\}\,m_h + \sum_{g<h,\,g\in\G} z_{gh}\, m_g\,,
\]
for some $z_{gh}\in\Q(q)$, and $\max(r_h) = \max(m_h) = h$ for
all $h\in\G$.
Finally, by Proposition~\ref{proptau}, $\max(u)=\max(\si(u))$
for all $u\in\U$, hence using Proposition~\ref{proport}, 
$\max(E_g)=\max(r_g)=g$ for all $g\in\G$.
\cqfd


\section{Canonical bases} \label{SECT3}

Kashiwara \cite{K} and Lusztig \cite{L} have introduced
independently and by different methods a canonical basis $\B$ of
$U_q(\nn)$.
Let $\Bs$ be the basis dual to $\B$ with respect to the scalar
product of \ref{ss1.2}. 
In this section we study the image of $\Bs$
in the embedding $\Phi : U_q(\nn) \longrightarrow \F$.

\subsection{}\label{newproof} 
The results of section~\ref{SECT2bis} give an easy alternative proof
of the existence of $\B$, as we shall now see.
For $g\in\G$ put
$
M_g:=\Phi(e_{\tau(g)}).
$
Inverting the formula of Lemma~\ref{PBW2mon}, we
get
\[
M_g=\sum_{h\in\G,\, h\ge g} \b_{gh}(q) \, E_h,
\]
where $\b_{gg}(q)=\a_{gg}(q)^{-1}=\k_g$.
Write
\[
\overline{E_g}=\sum_{h\in\G} a_{gh}(q)\,E_h,\qquad (g\in\G).
\]
Since $\overline{\U_\A}=\U_\A$ and $\{E_g\}$ is an $\A$-basis
of $\U_\A$, the coefficients $a_{gh}(q)$ belong to $\A$.
\begin{lemma}\label{involution}
$a_{gg}(q)=1$ for all $g\in\G$, and $a_{gh}(q)=0$ if $g>h$.
\end{lemma}
\proof
Clearly we have $\overline{M_g}=M_g$ for all $g\in\G$.
It follows that 
\[
a_{gh}(q)=\sum_{h\ge k\ge g;\,k\in\G} \a_{gk}(q^{-1}) \b_{kh}(q) \,.
\]
Hence $a_{gh}(q)=0$ if $g>h$, and 
$a_{gg}(q)=\a_{gg}(q^{-1})\b_{gg}(q)=\a_{gg}(q)\b_{gg}(q)=1$ 
by Lemma~\ref{PBW2mon}.
\cqfd

Let $L$ be the $\Z[q]$-lattice spanned by $\{E_g\}$.
It is well-known that Lemma~\ref{involution}
implies for all $g\in\G$ the existence of a unique $b_g\in L$ of the form
\begin{equation}\label{B2PBW}
b_g=E_g + \sum_{h\in\G,\, h> g} \gamma_{gh}(q) \, E_h,
\end{equation}
such that
$\gamma_{gh}(q)\in q\Z[q]$ and
\begin{equation}\label{barinv}
\overline{b_g}=b_g\,,
\end{equation}
(see for example \cite{L90} 7.10).
Clearly, $\{b_g\ |\ g\in\G\}$ is a $\Q(q)$-basis of $\U$,
a $\Z[q,q^{-1}]$-basis of $\U_\A$, and a $\Z[q]$-basis of $L$.

By Equations (\ref{SP1}), (\ref{SP2}) we have 
$(E_g\,,E_g)_{\{q=0\}} = 1$, hence by (\ref{B2PBW})
\begin{equation}\label{almostorth}
(b_g\,,b_h)_{\{q=0\}} = \delta_{gh}\,.
\end{equation}
It is easy to see that (up to sign) there is a unique 
$\A$-basis of $\U_\A$ satisfying (\ref{barinv})
and (\ref{almostorth}) (see \cite{L} 14.2).
Therefore, although the basis $\{E_g\}$ depends
on the choice of a total order on 
the set $\Pi$ of simple roots of~$\g$, the basis
$\{b_g\}$ is independent of this order.
This is the image under $\Phi$ of the canonical basis $\B$.

Equation~(\ref{B2PBW}) yields the next proposition,
which is needed for the proof of Theorem~\ref{th7}.
\begin{proposition}\label{triang}
For any total order on $\Pi$, the transition matrix
from $\{E_g\}$ to $\{b_g\}$ is unitriangular, if one arranges its 
rows and columns in lexicographic order.
\cqfd
\end{proposition}
For types $A$, $D$, $E$, a similar unitriangularity result was proved 
by Lusztig \cite{L90} for the PBW-type bases coming from the theory of Hall algebras.
As noted in \ref{PBWmain} the convex orderings of $\De^+$ coming from Lyndon
coverings are in general different from those coming from
Hall algebras. 
Hence, even for simply laced type, Proposition~\ref{triang} is 
different from Lusztig's result.

\subsection{}
Denote by $\{E^*_g\}$ the basis of $\U$ adjoint to
$\{E_g\}$, and by $\{b^*_g\}$ the basis adjoint to
$\{b_g\}$.
These are the images under $\Phi$ of the 
dual PBW-type basis and the dual canonical basis, respectively.

\begin{proposition}\label{charact}
The vector $b^*_g$ is characterized by the two following properties:

\noindent
{\rm (i)}\quad $b^*_g-E^*_g$ is a linear combination of vectors $E^*_h$
with coefficients in $q\Z[q]$;

\noindent
{\rm (ii)}\ \ The coefficients of the expansion of $b^*_g$ on the
basis of words are symmetric in $q$ and $q^{-1}$.
\end{proposition}
\proof
Clearly, $b^*_g$ satisfies (i). For $g, h \in \G$ we have
\[
(\si(b^*_g)\,,b_h)=(\tau(\overline{b^*_g})\,,b_h)
=q^{-N(|g|)}\,\{\overline{b^*_g}\,,b_h\}
=q^{-N(|g|)}\,\overline{(b^*_g\,,\overline{b_h})}
=q^{-N(|g|)}\,\delta_{gh}\,.
\]
Hence $\si(b^*_g)=q^{-N(|g|)}\,b^*_g$ and, by Lemma~\ref{symq},
$b^*_g$ satisfies (ii).
Now if $v$ is another element of $\U$ satisfying (i),
then 
$v=b^*_g+\sum_{h\not = g} \ga_h(q)b^*_h$ for  
some $\ga_h(q)\in q\Z[q]$.
If moreover $v$ satisfies (ii) then $\ga_h(q^{-1})=\ga_h(q)$,
hence $\ga_h(q)=0$ for all $h\in\G$, and $v=b^*_g$.
\cqfd

\begin{theorem}\label{th7}
We have $\max(b_g^*)=g$ for all $g\in\G$.
Moreover, the coefficient of the word $g$ in $b_g^*$ is equal to~$\k_g$.
\end{theorem}
\proof
We have
$E_g = \sum_{h\in\G} (E_g\,,b^*_h)\,b_h$.
By Proposition~\ref{triang},
$(E_g\,,b^*_h) = 0$ if 
$g>h$, and $(E_g\,,b^*_g)=1$.
Hence 
\begin{equation}\label{triang*}
b^*_h=E^*_h + \sum_{g\in\G;\,g< h} (E_g\,,b^*_h)\,E^*_g,
\end{equation}
and the results follow from 
\[
E^*_g=\sum_{k\le g} \b_{kg}(q) m_k = \k_gm_g + \sum_{k< g} \b_{kg}(q) m_k\,.
\]

\cqfd

Note that Theorem~\ref{th7} holds for any of the $r!$ different
total orders on the set $\Pi$ of simple roots of~$\g$.
The words $w$ such that $w=\max b^*$ for some total order
on $\Pi$ are similar to the extremal weights of the irreducible
character of a $\g$-module.


Finally, we note the following obvious consequence of
Proposition~\ref{triang} and Corollary~\ref{smallest}:
\begin{corollary}
For each $l\in\GL$ we have $E^*_l=b^*_l$, that is, the
root vector $E^*_l$ belongs to the dual canonical basis.
\cqfd
\end{corollary}

\subsection{}
We have the following important positivity property of $\Phi(\Bs)$.
\begin{theorem}\label{positive}
Assume that $\g$ is of type $A$, $D$, or $E$.
For all $g\in\G$, the coefficients $D_w(b^*_g)$
of the expansion of $b^*_g$ on the basis $\{w\in\M\}$
of $\F$ belong to $\N[q,q^{-1}]$.
\end{theorem}
\proof
Let $b=\Phi^{-1}(b_g)\in\B$ and $b^*=\Phi^{-1}(b^*_g)\in\Bs$.  
For $w=w[i_1,\ldots,i_k]$ with $|w|=|g|$ we have
\[
D_w(b^*_g) = (\Phi(e_{\tau(w)})\,,b^*_g)=(e_{\tau(w)}\,,b^*)\,.
\]
This is the coefficient of $b$ in the $\B$-expansion of
$e_{i_k}\cdots e_{i_1}$.
If $\g$ is of type $A, D, E$, Lusztig has shown
that these coefficients belong to $\N[q,q^{-1}]$ (\cite{L} 14.4.13).
\cqfd

In the non-simply laced case, the structure
constants of the multiplication on the canonical basis of
$U_q(\nn)$ need not be positive in general (see \cite{L} 14.4.14).
Nevertheless, the following conjecture is supported by extensive 
calculations.
\begin{conjecture}\label{conjpos}
For $\g$ of type $B, C, F_4$ or $G_2$ the coefficients $D_w(b^*_g)$
of the expansion of $b^*_g$ on the basis $\{w\in\M\}$
of $\F$ belong to $\N[q,q^{-1}]$.
\end{conjecture}     
%

Note that the conjecture can easily be proved in type $B_2$, using
the known fact that in this case all elements of $\Bs$ are 
$q$-commutative monomials in $8$ prime elements \cite{RZ,C}.
The images of these elements under $\Phi$ are
\[
w[1],\ \ w[2],\ \ w[1,2],\ \ w[2,1],\ \ [2]\,w[1,1,2],\ \ [2]\,w[2,1,1],
\ \ w[1,2,1],\ \ [2]\,w[2,1,1,2].
\]
Clearly, all $q$-shuffle monomials in these words have
positive coefficients.

Following Lusztig \cite{L3}, define the variety $N_{\ge 0}$ of 
totally nonnegative elements in $N$ as the monoid generated by
the $x_i(t) \ (1\le i\le r,\ t\in\R_{\ge 0})$.
It follows from Proposition~\ref{phi_inverse} that if 
all the coefficients of the $\M$-expansion of $\varphi(f)$
belong to $\R_{\ge 0}$, then $f$ takes nonnegative values on $N_{\ge 0}$.
In particular, for simply-laced type, we recover by means of 
Theorem~\ref{positive} the known fact that the functions of $\C[N]$
obtained by specializing at $q=1$ the elements of $\Bs$ 
take nonnegative values on $N_{\ge 0}$.
This property is also true for non simply-laced type.
Note however that one can easily find examples of functions $f$
which are nonnegative on $N_{\ge 0}$ while $\varphi(f)$ has
some negative coefficients.

In fact, Proposition~\ref{varphi} and \ref{phi_inverse} show
that all the coefficients of $\varphi(f)$ are nonnegative
if and only if all the coefficients of the polynomial function
$(t_1,\ldots ,t_k) \mapsto f(x_{i_1}(t_1)\cdots x_{i_k}(t_k))$
are nonnegative for any sequence $(i_1,\ldots , i_k)$.
Thus Conjecture~\ref{conjpos} at $q=1$ can be formulated in the
following way:
let $G$ be of non-simply laced type and
let $f_b\in\C[N]$ denote the specialization at $q=1$ of $b\in\Bs$;
then $f_b(x_{i_1}(t_1)\ldots x_{i_k}(t_k))\in \N[t_1,\ldots ,t_k]$
for all sequences $(i_1,\ldots ,i_k)$.
This has been proved by Berenstein and Zelevinsky in the case
where $f_b$ is a generalized minor (\cite{BZ2} Th. 5.8).

\subsection{}\label{indecomp}
Recall from \ref{special2} that $c_w^{-1}e_w$ denotes a monomial
in the divided powers of the Chevalley generators.
In this section, we shall say for short that $c_w^{-1}e_w$
is `a monomial'. 
\begin{lemma}\label{cof}
{\rm (i)} $\Phi(\Bs)$ is an $\A$-basis of $\U^*_\A$. 

\noindent
{\rm (ii)} Suppose that the monomial $c_w^{-1}e_{\tau(w)}$ belongs to
$\B$, and let $b_g=\Phi(c_w^{-1}e_{\tau(w)})$.
Then $w$ occurs in the $\M$-expansion of $b^*_g$ with 
coefficient $c_w w$, and for all $h\not = g$, 
$w\not \in \SS(b^*_h)$.
\end{lemma}
\proof
(i) It is known that $\B$ is an $\A$-basis of $U_\A$.
It follows that $\B^*$ is an $\A$-basis of $U^*_\A$
and $\Phi(\Bs)$ is an $\A$-basis of $\U^*_\A=\Phi(U^*_\A)$
(see Lemma~\ref{Abasis}).

(ii) We have 
$\de_{gh}=(b_g\,,b^*_h)=c_w^{-1}(e_{\tau(w)}\,,b^*_h),$
and this is the coefficient of $w$ in the $\M$-expansion of 
$c_w^{-1}\,b^*_h$.
\cqfd

For $\nu\in Q^+$ set
$
\B_\nu=\{b\in\B\ |\ |b|=\nu\},
$
$
\Bs_\nu=\{b^*\in\B^*\ |\ |b^*|=\nu\},
$
and
\[
\F_\N(\nu)=\bigoplus_{w\in\M,\,|w|=\nu}\N[q,q^{-1}]\,c_w w\,,
\qquad
\U^*_\N(\nu)=\bigoplus_{g\in\G,\,|g|=\nu}\N[q,q^{-1}]\,b^*_g\,.
\]
We shall assume until the end of section~\ref{indecomp} that 
$\U^*_\N(\nu) \subset \F_\N(\nu)$. 
(By Theorem~\ref{positive} this holds for all $\nu$ in type $A$, $D$, $E$.)
\begin{proposition}\label{minimal}
The following statements are equivalent:

\noindent
{\rm (i)}\ \  all elements of $\B_\nu$ are monomials;

\noindent
{\rm (ii)}\ \ for every $b^*_g\in\Phi(\B^*_\nu)$ there exists $w\in\M$
such that $w$ occurs in the $\M$-expansion of $b^*_g$ with coefficient
$c_w$ and $w$ does not occur in any other $b^*_h$; 

\noindent
{\rm (iii)}\ \  $\U^*_\N(\nu) = \U \cap \F_\N(\nu)$.
\end{proposition}
\proof
By Lemma~\ref{cof}~(ii), we have that (i) implies (ii).
Conversely, if (ii) holds, the monomials $c_w^{-1}e_{\tau(w)}$
associated with each $g$ form a family of vectors adjoint
to $\B^*_\nu$, and (i) follows.

Suppose that (ii) holds, and let 
\[
v=\sum_{w\in\M,\,|w|=\nu} \si_{vw}\,c_w w
=\sum_{g\in\G,\,|g|=\nu} \tau_{vg}\,b^*_g
\]
be an element of $\U \cap \F_\N(\nu)$.
Then for each $g$ there exists a $w$ such that
$\tau_{vg}=\si_{vw}\in\N[q,q^{-1}]$, and (iii) holds.

Suppose that (ii) does not hold. 
Then there exists $g_0$ such that every $w$ occuring in the $\M$-expansion
of $b^*_{g_0}$ occurs also in the expansion of some $b^*_{g_w}$ for some 
$g_w\not = g_0$.
Write 
\[
b^*_g=\sum_{w\in\M}\varphi_{gw}\,w\,,\qquad (g\in\G, \ |g|=\nu).
\]
Then, denoting by $\de$ the l.c.m. of all $\varphi_{g_w w}$ for
$w\in\SS(b^*_{g_0})$, and setting 
$\de_w:=\de\,\varphi_{g_0w}\,\varphi_{g_ww}^{-1}$
we have that
\[
-\de\,b^*_{g_0}+\sum_{w\in\SS(b^*_{g_0})} \de_w\,b^*_{g_w}
\]
belongs to $\U \cap \F_\N(\nu)$ but not to $\U^*_\N(\nu)$.
\cqfd

\begin{definition}
We say that $v\in\U \cap \F_\N(\nu)$ is indecomposable
if {\rm (a)} there exists no decomposition
$v=v_1+v_2$ with nonzero $v_1, v_2 \in\U \cap \F_\N(\nu)$,
and {\rm (b)} there exists $w\in\SS(v)$ whose coefficient
is equal to $c_w$. 
\end{definition}
Concretely, what Proposition~\ref{minimal} means is that
$\B_\nu$ consists only of monomials if and only if $\Phi(\B^*_\nu)$
consists of all the indecomposable elements of $\U \cap \F_\N(\nu)$. 

It is well known that all elements of $\B$ are monomials when
$\g$ is of type $A_2$, hence in this case $\Phi(\Bs)$ is precisely
the set of all indecomposable elements of $\U \cap \F_\N(\nu)$. 
In general, there are indecomposable elements which do not
belong to $\Phi(\B^*)$. 
It may also happen that some elements of $\Phi(\B^*)$ are not indecomposable.
It seems to be an interesting problem to understand which elements
of $\Phi(\B^*)$ are indecomposable. 

\subsection{}\label{algo1}
We now describe an algorithm to compute the basis $\{b^*_g\}$.
All calculations take place in the $q$-shuffle algebra $(\F,\shuf)$
and all vectors are expressed on the basis $\M$ of words.
We fix an arbitrary total order on~$\Pi$. 

\subsubsection{}
The first step is to calculate the set $\GL$ of good Lyndon words.
For this we use \ref{LyndCover}. 

\subsubsection{}
For each $l\in\GL$ we calculate $r_l$ as an iterated
$q$-bracket given by the co-standard factorization of $l$.
Then we obtain $E^*_l$ by an appropriate normalization
of $r_l$.
Namely, we have 
\begin{equation}\label{kappal}
\k_l E^*_l = {(-1)^{\ell(l)-1}\over q^{N(|l|)} (E_l,E_l)}\, r_l\,,
\end{equation}
where $(E_l,E_l)$ is given by (\ref{SP2}).
It remains to calculate $\k_l$.
By Theorem~\ref{th7}, we know that the coefficient of $l$
in $E^*_l$ is equal to $\k_l$.
Hence the coefficient of $l$ in (\ref{kappal}) is equal to
$\k_l^2$, and to get $E^*_l$ we just need to divide 
(\ref{kappal}) by the square root of its coefficient of $l$.

\subsubsection{}
Let us fix a weight $\nu\in Q^+$.
By Proposition~\ref{factgood}, we can easily calculate the 
ordered list $\{g_1<\ldots < g_s\}$ of all good words of
weight $\nu$.
Note that for a good word $g=l_1^{a_1}\cdots l_k^{a_k}$
with $l_1>\cdots >l_k\in\GL$ we have
\[
E^*_g=q^{c_g}\,
(E^*_{l_k})^{*a_k}\shuf\cdots \shuf (E^*_{l_1})^{*a_1}\,
\]
where $c_g=\sum_{i=1}^k{a_i\choose 2}d_{l_i}$
(this follows easily from (\ref{SP1}) (\ref{SP2})).
So we can compute $E^*_{g_1},\ldots , E^*_{g_s}$.
By (\ref{triang*}), we have $b^*_{g_1}=E^*_{g_1}$.
Suppose that for some $t\le s$ we have calculated 
$b^*_{g_1},\ldots ,b^*_{g_{t-1}}$.
If all the coefficients of the expansion of $E^*_{g_t}$
on the basis of words are symmetric in $q$ and $q^{-1}$
then $b^*_{g_t}=E^*_{g_t}$. 
Otherwise let $g_j$ be the largest good word occuring
in $E^*_{g_t}$ with a non symmetric coefficient $\a\in\Z[q,q^{-1}]$.
We know that the coefficient of $g_j$ in $b^*_{g_j}$ is
$\k_{g_j}$, which is symmetric in $q$ and $q^{-1}$.
The existence of $b^*_{g_t}$ implies that there exists $\ga\in q\Z[q]$
such that the coefficient $\rho$ of $g_j$ in
$E^*_{g_t}-\ga b^*_{g_j}$ is symmetric in $q$ and $q^{-1}$.
Moreover if there were other coefficients $\ga'$ and $\rho'$
satisfying the same properties, we would have 
$\a=\k_{g_j}\ga + \rho = \k_{g_j}\ga' + \rho'$, hence
$\k_{g_j}(\ga-\ga')=\rho'-\rho$, with $\k_{g_j}$, $\rho'-\rho$
symmetric in $q$ and $q^{-1}$ and $\ga-\ga'\in q\Z[q]$.
This forces $\rho'-\rho=\ga-\ga'=0$, therefore $\ga$ is uniquely
determined.
If now all the coefficients of the expansion of $E^*_{g_t}-\ga b^*_{g_j}$
on the basis of words are symmetric in $q$ and $q^{-1}$
then $b^*_{g_t}=E^*_{g_t}-\ga b^*_{g_j}$, otherwise
we apply the same procedure as above to $E^*_{g_t}-\ga b^*_{g_j}$. 
After a finite number of steps we will obtain $b^*_{g_t}$. 

\subsubsection{}\label{EXG2}
Let us demonstrate the algorithm on an example.
We choose $\g$ of type $G_2$.
Then 
\[
\Delta^+ = \{\a_1,\,\a_2,\,\a_1+\a_2,\,2\a_1+\a_2,\,3\a_1+\a_2,\,3\a_1+2\a_2\}\,.
\]
Let us decide that $w_1<w_2$.
Then, the procedure of \ref{LyndCover} gives immediately
\[
\GL=\{w[1], w[1,1,1,2], w[1,1,2], w[1,1,2,1,2], w[1,2], w[2]\}\,.
\]
Let us calculate for example the vector $E^*_{w[1,1,2,1,2]}$.
We have
\[
r_{w[1,1,2,1,2]}=r_{w[1,1,2]}\shuf r_{w[1,2]} - q\,r_{w[1,2]}\shuf
r_{w[1,1,2]}\,,
\]
and by induction, we may assume that $r_{w[1,1,2]}$
and 
$r_{w[1,2]}$
are already known.
Rescaling as indicated above we get
\[
\k_{w[1,1,2,1,2]} E^*_{w[1,1,2,1,2]} =
[2]^2_1[3]^2_1\,w[1,1,2,1,2] + [2]^2_1[3]^2_1[2]_3\,w[1,1,1,2,2]\,,
\]
hence $\k_{w[1,1,2,1,2]}=[2]_1[3]_1$ and 
\[
E^*_{w[1,1,2,1,2]} = [2]_1[3]_1\,w[1,1,2,1,2] +
[2]_1[3]_1[2]_3\,w[1,1,1,2,2]\,.
\]
The other root vectors are calculated similarly and one finds
\[
E^*_{w[1,1,1,2]}=[2]_1[3]_1\,w[1,1,1,2],\quad
E^*_{w[1,1,2]}=[2]_1\,w[1,1,2],\quad
E^*_{w[1,2]}=w[1,2]\,.
\]
Let us calculate the dual canonical basis of the 
weight space corresponding to the highest root 
$\b=3\a_1+2\a_2$.
The list of good words of weight $\b$ in increasing order is
\begin{eqnarray*}
\G_\b&=&\{w[1,1,2,1,2],w[1,2,1,1,2],w[1,2,1,2,1],\\
&&\ \ w[2,1,1,1,2],w[2,1,1,2,1],w[2,1,2,1,1],w[2,2,1,1,1]\}\,.
\end{eqnarray*}
We have already calculated $b^*_{w[1,1,2,1,2]}=E^*_{w[1,1,2,1,2]}$.
Next, we have
\begin{eqnarray*}
E^*_{w[1,2,1,1,2]}&=&E^*_{w[1,1,2]}\shuf E^*_{w[1,2]}\\
&=&[2]_1\,w[1,2,1,1,2]+q[2]_1[3]_1\,w[1,1,2,1,2]
+q[2]_1[3]_1[2]_3\,w[1,1,1,2,2]\,. 
\end{eqnarray*}
Hence 
$
b^*_{w[1,2,1,1,2]}=E^*_{w[1,2,1,1,2]}-q\,b^*_{w[1,1,2,1,2]}
=[2]_1\,w[1,2,1,1,2]\,.
$
Next, we have
\begin{eqnarray*}
E^*_{w[1,2,1,2,1]}&=&qE^*_{w[1]}\shuf E^*_{w[1,2]}\shuf E^*_{w[1,2]}\\
&=&[2]_1\,w[1,2,1,2,1]+q^2[2]_1^2\,w[1,2,1,1,2]+[2]_1[2]_3\,w[1,1,2,2,1]\\
&&\ \ +\,([2]_1+q^4[2]_1[2]_3)\,w[1,1,2,1,2]
+q^4[2]_1[2]_3[3]_1\,w[1,1,1,2,2]\,,
\end{eqnarray*}
hence
\begin{eqnarray*}
b^*_{w[1,2,1,2,1]}&=&E^*_{w[1,2,1,2,1]}-q^2\,[2]_1\,b^*_{w[1,2,1,1,2]}
-q^4\,b^*_{w[1,1,2,1,2]}\\
&=&[2]_1\,w[1,2,1,2,1]+[2]_1\,w[1,1,2,1,2]+[2]_1[2]_3\,w[1,1,2,2,1]\,.
\end{eqnarray*}
In the same way one calculates
\begin{eqnarray*}
b^*_{w[2,1,1,1,2]}&=&[2]_1[3]_1\,w[2,1,1,1,2],\\
b^*_{w[2,1,1,2,1]}&=&[2]_1\,w[2,1,1,2,1],\\ 
b^*_{w[2,1,2,1,1]}&=&[2]_1\,w[2,1,2,1,1]+[2]_1\,w[1,2,1,2,1]+
[2]_1[2]_3\,[1,2,2,1,1],\\ 
b^*_{w[2,2,1,1,1]}&=&[2]_1[2]_3[3]_1\,w[2,2,1,1,1]+[2]_1[3]_1w[2,1,2,1,1]\,.
\end{eqnarray*}


\section{Type $A$ and $q$-characters of affine Hecke algebras} \label{SECT4}

\subsection{} Let $t\in\C^*$ be of infinite multiplicative order.
Let $H_m=H_m(t)$ be the algebra over $\C$ generated by 
invertible elements $T_1,\ldots , T_{m-1},y_1,\ldots , y_m$ 
subject to the following relations:
\begin{eqnarray*}
&&T_iT_{i+1}T_i=T_{i+1}T_iT_{i+1},\hskip 1.1cm (1\le i\le m-2),\label{EQ_T1}\\
&&T_iT_j=T_jT_i,\hskip 2.9cm (\vert i-j\vert>1),\label{EQ_T2}\\
&&(T_i-t)(T_i+1)=0, \hskip 1.6cm(1\le i\le m-1),\label{EQ_T3}\\
&&y_iy_j=y_jy_i, \hskip 3.1cm (1\le i,j\le m),\label{EQ_YY}\\
&&y_jT_i=T_iy_j, \hskip 3cm (j \not= i,i+1),\\
&&T_iy_iT_i=t\,y_{i+1}, \hskip 2.5cm (1\le i\le m-1).\label{EQ_Y}
\end{eqnarray*}
This is the Bernstein presentation of the affine Hecke algebra of
$GL(m)$. 
\subsection{}
Let $M$ be a finite-dimensional $H_m$-module.
Since the elements $y_i$ are pairwise commutative,
$M$ decomposes as a sum of generalized eigenspaces
\[
M=\bigoplus_\gamma M[\gamma]\,,
\]
where for $\gamma\in\C^m$, we put 
\[
M[\gamma]=\{m\in M\ |\ \mbox{\ for all i,\ }(y_i-\gamma_i)^{n_i} m = 0 
\mbox{\ for some\ } n_i\in\N^*\}\,.
\]
The $\ga$ such that $M[\ga]\not = 0$ are called the weights
of $M$. 
We will say that $M$ is integral if all its weights are of
the form $\ga = (t^{i_1},\ldots , t^{i_m})$ for some
$i_1,\ldots , i_m\in\Z$.
In that case we shall write $M[i_1, \ldots , i_r]$ in
place of $M[\ga]$.

\subsection{}
Let $\CC_{m,r}$ denote the category of integral $H_m$-module
with weights $(t^{i_1},\ldots , t^{i_m})$ such that
$1\le i_k\le r$ for all $k=1,\ldots ,m$.
The character of $M$ is defined by
\[
\ch M = \sum_{1\le i_1,\ldots ,i_m\le r} \dim M[i_1,\ldots ,i_m]\,w[i_1,\ldots ,i_m]\,.
\] 
This is an element of $\F_\C$ (see \ref{special1}).

\subsection{}\label{charshuffle}
Let $m=m_1+m_2$. The parabolic subalgebra $H_{m_1,m_2}$ of $H_m$
generated by
\[ 
T_1,\ldots ,T_{m_1-1},T_{m_1+1},\ldots ,T_{m-1},y_1,\ldots , y_m,
\]
is isomorphic to $H_{m_1}\otimes H_{m_2}$.
Let $M_1$ and $M_2$ be a $H_{m_1}$-module and a $H_{m_2}$-module,
respectively.
The induction product $M_1\odot M_2$ is the $H_m$-module defined
by 
\[
M_1\odot M_2 = \Ind_{H_{m_1,m_2}}^{H_m} M_1\otimes M_2\,.
\]
If $M_1$ and $M_2$ are objects of $\CC_{m_1,r}$ and $\CC_{m_2,r}$
then $M_1\odot M_2$ is an object of $\CC_{m,r}$ and we have \cite{GV}
\begin{equation}
\ch M_1\odot M_2 = \ch M_1 \shuffle \ch M_2\,,
\end{equation}
where $\shuffle$ is the classical shuffle product. 
This follows from a Mackey-type theorem for $H_m$.

Let $\RR = \bigoplus_{m\in\N} \RR_{m,r}$, where $\RR_{m,r}$
is the complexified Grothendieck group of $\CC_{m,r}$
(by convention, we put $\RR_{0,r}=\C$).
The class in $\RR$ of a module $M$ is denoted by $[M]$.
The operation $\odot$ induces in $\RR$ a multiplication $\times$
that makes it into a $\C$-algebra.
Note that $\times$ is commutative: although 
$M_1\odot M_2$ is in general not isomorphic to $M_2\odot M_1$,
their classes in $\RR$ coincide.
Then 
\[
\ch : (\RR,\times) \longrightarrow (\F_\C,\shuffle)
\]
is a ring homomorphism.

\subsection{}\label{polynom}
For $1\le i\le j\le r$, let $M_{[i,j]}$ be the $1$-dimensional
$H_{j-i+1}$-module on which the $T_k$'s act by multiplication
by $t$, and the $y_k$'s by multiplication by $t^{k+i-1}$.
It is known that $\RR$ is the polynomial ring over $\C$
in the variables $[M_{[i,j]}]\ (1\le i\le j\le r)$ \cite{Z}. 
Now,
$
\ch M_{[i,j]} = w[i,\ldots,j]\,.
$
Therefore, $\ch \RR$ is the subring of $(\F_\C,\shuffle)$
generated by the words $w[i,\ldots,j] \ (1\le i\le j\le r)$.

\subsection{}
A multi-segment $\mm$ is a list of segments
$
\mm=([i_1,j_1],\ldots ,[i_k,j_k])
$
written in increasing order with respect to the
following total order on segments:
\[
[i,j]<[k,l] \quad\Longleftrightarrow\quad
(i<k \mbox{ or } (i=k \mbox{ and } j<l))\,.
\]
Following Zelevinsky\cite{Z}, to $\mm$ we associate a standard induced module
\[
M_\mm = M_{[i_1,j_1]}\odot \cdots \odot M_{[i_k,j_k]}
\] 
and a simple module $L_\mm$ (see for example \cite{Rog} or \cite{LNT}).

Note that the words $w[i,\ldots,j] \ (1\le i\le j\le r)$ are
the good Lyndon words for the root system $A_r$
corresponding to the natural order
$w_1<\cdots < w_r$, and the multi-segments $\mm$ are in
one-to-one correspondence with the good words $g$
by 
\begin{equation}\label{mg}
\mm=([i_1,j_1],\ldots ,[i_k,j_k])\ \longleftrightarrow \ 
g=w[i_k,\ldots ,j_k,i_{k-1},\ldots , j_{k-1},\ldots ,i_1,\ldots ,j_1]
\,.
\end{equation}

\subsection{}
Let $\g=\Sl_{r+1}$ be the Lie algebra of type $A_r$, and let 
$U_q(\nn)$ be the corresponding quantum algebra.
Choose the convex ordering $\b_1<\cdots <\b_n$ of $\Delta^+$ associated
with the reduced decomposition
\[
w_0=s_1s_2s_1s_3s_2s_1\cdots s_{r-1}\cdots s_2s_1\,.
\]
It is easy to check that this is the same as the convex
ordering coming from the good Lyndon words above,
namely 
\[
\a_i+\cdots +\a_j < \a_k+\cdots +\a_l
\quad\Longleftrightarrow\quad
w[i,\ldots ,j] < w[k,\ldots ,l]\,.
\]
The PBW-type basis of $U_q(\nn)$ associated with this choice is thus conveniently
labelled by multi-segments $\mm=\sum_{1\le i,\le j\le r} m_{ij}\,[i,j]$
where $m_{ij}$ denotes the multiplicity of the segment $[i,j]$.
We shall write
\[
E_\mm:=\prod_{1\le i \le j\le r} E(\a_i+\ldots+\a_j)^{(m_{ij})}
\]
where the product is taken in the order given by the convex ordering
above.
We denote accordingly $\{E^*_\mm\}$ and $\{b^*_\mm\}$ the dual PBW-basis
and dual canonical basis respectively.
We have $\Phi(E^*_\mm)=E^*_g$ and $\Phi(b^*_\mm)=b^*_g$ where the 
correspondence between multisegments $\mm$ and good words $g$ is
given by (\ref{mg}).
Moreover, it is easy to check that 
$\Phi(E^*_{[i,j]})=w[i,\ldots ,j]$.

\subsection{}\label{specializ} 
Recall the setup of \ref{special1}. 
Let $\underline{E}^*_\mm\in\C[N]$ and $\underline{b}^*_\mm\in\C[N]$ denote
the specializations of $E^*_\mm$ and~$b^*_\mm$ at $q=1$.
Then $\underline{E}^*_{[i,j]}\ (1\le i\le j\le r)$ is just 
the coordinate function $t_{i,j+1}$ mapping a matrix $g$ to its entry $g_{i,j+1}$.
It follows from \ref{polynom} that $\C[N]$ is isomorphic as an algebra 
to $(\RR,\times)$. 
Let  $\theta : \C[N] \longrightarrow \RR$ denote this isomorphism.
By a dual version of Ariki's theorem (\cite{A}, see also \cite{LNT}), 
we have more precisely
\begin{equation}
\theta(\underline{E}^*_\mm) = [M_\mm],\qquad 
\theta(\underline{b}^*_\mm) = [L_\mm]\,.
\end{equation}
Consider the diagram
\[
\begin{array}{ccccc}
U_\A^*&&\stackrel{\Phi}{\longrightarrow} && \F_\A\\[2mm]
{\downarrow}  &&&                &\downarrow\\[2mm]
\C[N] &  &\stackrel{\varphi}{\longrightarrow}& & \F_\C \\[2mm]
{\qquad \theta}&\searrow&  &\nearrow & {\ch\qquad} \\[2mm]
&&\RR&&
\end{array}
\]
where the two vertical arrows denote specialization $q\mapsto 1$.
For $\mm=([i_1,j_1],\ldots ,[i_k,j_k])$, we have 
\[
\varphi(\underline{E}^*_\mm))= 
w[i_1,\ldots ,j_1]\shuffle \cdots \shuffle w[i_k,\ldots,j_k]
= \ch M_\mm\,,
\]
hence the diagram is commutative. 
Therefore, for all multi-segments $\mm$,
\[
\ch M_\mm = {\Phi(E^*_\mm)}_{\{q=1\}},\qquad \ch L_\mm = {\Phi(b^*_\mm)}_{\{q=1\}}\,.
\]
In other words
\begin{theorem}\label{qcharA}
$\Phi(E^*_\mm)$ is a $q$-analogue of the character
of the standard induced module $M_\mm$ and
$\Phi(b^*_\mm)$ is a $q$-analogue of the character
of the simple module $L_\mm$.
\cqfd
\end{theorem}
It follows immediately from the definition of $\shuf$ that
the $q$-analogues $\Phi(E^*_\mm)$
of the characters of the standard modules have nonnegative
coefficients.  
By Theorem~\ref{positive}, this is also true for the $q$-analogues 
$\Phi(b^*_\mm)$ of the characters of the irreducible modules.

\subsection{} 
As an application of Theorem~\ref{qcharA}, we can calculate
explicitly a family of vectors $b^*_{g(\l/\m;s)}$
labelled by skew Young diagrams $\l/\m$ and an integer $s$.
  
Given two partitions $\l=(\l_1\ge\l_2\ge\cdots\ge\l_j>0)$
and  $\m=(\m_1\ge\m_2\ge\cdots\ge\m_k>0)$ such that $j\ge k$ and
$\l_i\ge\m_i$ for $i=1,\ldots ,k$, we denote by $\l/\m$ the 
skew Young diagram obtained by removing the Young diagram
of $\m$ from that of $\l$.
We define the content of the cell on row $i$ and column $j$
of $\l/\m$ to be $c=j-i$.
Let $m=\sum_{i=1}^j\l_i-\sum_{i=1}^k\m_i$.
A standard Young tableau $T$ of shape $\l/\m$
is a filling of the cells of $\l/\m$
by the integers $1,2,\ldots , m$, increasing
on rows and columns.
To $T$ and an integer $s$ we associate the word 
$w[T,s]:=w[c_1+s,\ldots ,c_m+s]$,
where $c_i$ denotes the content of the cell numbered $i$
in $T$.

\begin{figure}
\begin{center}
\leavevmode
\epsfxsize =6.8cm
\epsffile{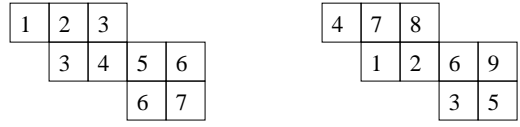}
\end{center}
\caption{\small \label{Fig0}{\it The skew Young diagram $\l/\m$ with
$\l=(5,5,3)$ and $\m=(3,1)$, filled with
its contents shifted by $s=3$, and a standard Young tableau $T$ of shape
$\l/\m$ with $w[T,3]=w[3,4,6,1,7,5,2,3,6]$.}}
\end{figure}

These definitions are illustrated in Figure~\ref{Fig0}.
Assume that $\l$ is such that $1-\l_1+r\ge j$ and let $s\in
[j,1-\l_1+r]$.
To $(\l/\m;s)$ we associate the good word
\[
g(\l/\m;s)
=w[\m_1+s,\ldots,\l_1-1+s,\m_2-1+s,\ldots,\l_2-2+s,\ldots,\m_j-j+1+s,\ldots,\l_j-j+s].
\]
This is the word obtained by reading the rows of $\l/\m$ from left to
right and bottom to top, the cells being filled by the contents
shifted by $s$.
(We assume that $\m$ is made into a sequence of length $j$ by
appending a tail of $j-k$ digits $0$.) 
\begin{corollary}\label{charskew}
We have 
\begin{equation}\label{tabl0}
b^*_{g(\l/\m;s)} = \sum_T w[T,s]\,,
\end{equation}
where $T$ runs through the set of all standard Young
tableaux of shape $\l/\m$.
\end{corollary}
\proof
To each choice of $\l/\m$ and $s$ as above corresponds an irreducible
$H_m$-module $L_{\mm(\l/\m;s)}$ on which the generators $y_1,\ldots , y_m$ act
semi-simply. 
The multi-segment $\mm(\l/\m;s)$ is obtained from the good word
$g(\l/\m;s)$ by the correspondence (\ref{mg}).
The character of $L_{\mm(\l/\m;s)}$ is known to be given
by the right-hand side of (\ref{tabl0}).
Moreover, the generalized eigenspaces of $L_{\mm(\l/\m;s)}$
are all $1$-dimensional.
Hence the $q$-character of $L_{\mm(\l/\m;s)}$ coincides
with its ordinary character, and the result follows from
Theorem~\ref{qcharA}.
\cqfd

Corollary~\ref{charskew} may also be proved directly (\ie without
using the representation theory of $H_m$) by arguing as in
Proposition~\ref{chartabshif} and Proposition~\ref{chartabshif2} below.


\section{Type $B$ and $q$-characters of affine Hecke-Clifford
  superalgebras} 
\label{SECT5}

\subsection{} Let us take $\g$ of type $B_r$.
We choose the following numbering of the simple roots
\begin{center}
\leavevmode
\epsfxsize =6.8cm
\epsffile{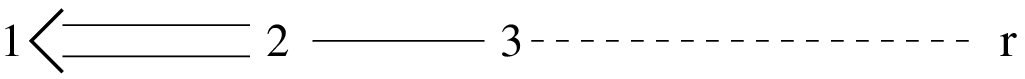}
\end{center}
and the standard ordering $w[1]<w[2]<\cdots <w[r]$.
The set of good Lyndon words is calculated using \ref{LyndCover},
and we find
\begin{equation}\label{GLB}
\GL=\{w[i,\ldots,j],\ 1\le i\le j\le r\}
\cup
\{w[1,\ldots,j,1,\ldots,k],\ 1\le j<k\le r\}.
\end{equation}
Here are some simple examples of vectors $b^*_g$. 
\begin{lemma}\label{examplesB}
\begin{eqnarray*}
b^*_{w[i,i+1]}&=&w[i,i+1],\hskip 5.4cm (1\le i\le r-1),\\
b^*_{w[i+1,i]}&=&w[i+1,i],\hskip 5.4cm (1\le i\le r-1),\\
b^*_{w[1,1,2]}&=&[2]_1\,w[1,1,2],\\
b^*_{w[1,2,1]}&=&w[1,2,1],\\
b^*_{w[2,1,1]}&=&[2]_1\,w[2,1,1],\\
b^*_{w[i,i+1,i]}&=&w[i,i+1,i]+[2]_2\,\,w[i,i,i+1],\hskip 2cm 
(2\le i\le r-1),\\
b^*_{w[i+1,i,i]}&=&[2]_2\,\,w[i+1,i,i]+w[i,i+1,i],\hskip 2cm 
(2\le i\le r-1),\\
b^*_{w[1,1,2,1]}&=&[2]_1\,w[1,1,2,1]+[3]_1[2]_1\,w[1,1,1,2],\\
b^*_{w[1,2,1,1]}&=&[2]_1\,w[1,2,1,1]+[2]_1\,w[1,1,2,1],\\
b^*_{w[2,1,1,1]}&=&[3]_1[2]_1\,w[2,1,1,1]+[2]_1\,w[1,2,1,1].\\
\end{eqnarray*}
\end{lemma}
\proof
These are straightforward calculations using the algorithm of \ref{algo1}. 
\cqfd

\subsection{}
\begin{figure}
\begin{center}
\leavevmode
\epsfxsize =6.8cm
\epsffile{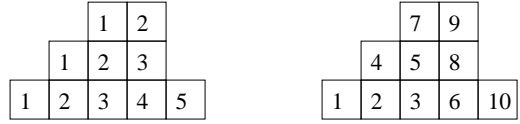}
\end{center}
\caption{\small \label{Fig1}{\it The shifted Young diagram of $\l=(5,3,2)$ filled with
its contents, and a standard shifted Young tableau $T$ of shape
$\l$ with $w[T]=w[1,2,3,1,2,4,1,3,2,5]$.}}
\end{figure}

Let $\l=(\l_1>\l_2>\cdots >\l_k>0)$ be a strict partition
with $\l_1\le r$.
We set
\[
\nu_\l=\sum_{i=1}^k\a_1+\cdots +\a_{\l_i}\in Q^+\,.
\]
We represent $\l$ graphically by a shifted Young diagram.
We define the content of the cell on row $i$ and column $j$
of a shifted Young diagram to be $c=j-i+1$.
A standard shifted Young tableau $T$ of shape $\l$
is a filling of the cells of the shifted diagram of $\l$
by the integers $1,2,\ldots , m=\sum_i \l_i$, increasing
on rows and columns.
To $T$ we associate the word $w[T]:=w[c_1,\ldots ,c_m]$,
where $c_i$ denotes the content of the cell numbered $i$
in $T$.
These definitions are illustrated in Figure~\ref{Fig1}.
Finally, we associate to $\l$ the good word
$
g(\l):=w[1,\ldots ,\l_1,1,\ldots ,\l_2,\ldots,1,\ldots,\l_k]\,.
$
\begin{proposition}\label{chartabshif}
We have 
\begin{equation}\label{tabl}
b^*_{g(\l)} = \sum_T w[T]\,,
\end{equation}
where $T$ runs through the set of all standard shifted Young
tableaux of shape $\l$.
\end{proposition}
\proof
Let $\L\in\h^*$.
Introduce the adjoint action twisted by $\L$ from $U_q(\g)$
to $\End U_q(\nn)$.
It is defined by
\begin{eqnarray*}
\Ad_{\L}(f_i)(x)&=&{1\over q^{d_i}-q^{-d_i}}
\left(q^{(\L\,,\,\a_i)}xe_i-q^{-(\L\,,\,\a_i)+(|x|\,,\,\a_i)}e_ix\right)\,,\\
\Ad_{\L}(e_i)(x)&=&e'_i(x)\,,
\end{eqnarray*}
for a homogeneous $x\in U_q(\nn)$.
It is well-known that $\Ad_\L$ endows $U_q(\nn)$ with the 
structure of a dual Verma module $M(\L)^*$ with highest weight $\L$.
Moreover the dual canonical basis $\Bs_\L$ of the irreducible
submodule $V(\L)$ generated by the highest weight vector
of $M(\L)$ becomes in this realization a subset of the 
dual canonical basis $\Bs$ of $U_q(\nn)$.
Let $\L_1$ be the first fundamental weight, so that $V(\L_1)$
is the spin representation.
This is a minuscule representation of dimension $2^r$ 
for which the canonical basis
and the dual canonical basis coincide and are given by  
\[
\B_{\L_1}=\B^*_{\L_1}=\{\Ad_{\L_1}(f_{i_1}\cdots f_{i_k})(1_{U_q(\nn)}) 
\ |\ 1\le i_1,\ldots,i_k \le r,\ 1\le k \le r(r+1)/2\} \setminus \{0\}\,.
\] 
We are going to prove that
 $\{b^*_{g(\l)}\}=\Phi(\Bs_{\L_1})=\{S_\l\}$,
where $S_\l$ denotes the tableau sum of (\ref{tabl}).

It is well known that $\Bs_{\L_1}$ has a natural indexation 
by the strict partitions $\l$ as above, namely we write 
$b^*_\l$ for the unique element of $\Bs_{\L_1}$ of weight
$\L_1-\nu_\l$.
Then we have 
\[
\Ad_{\L_1}(e_i)(b^*_\l) = 
\left\{
\matrix{
b^*_{(\l_1,\ldots ,\l_j-1,\ldots ,\l_k)} &\mbox{\rm if $\l_j= i$ and
  $\l_{j+1}\not = i-1$,}\cr
0\hfill &\mbox{otherwise.}\hfill\hfill
}
\right.
\] 
Using Theorem~\ref{th2}, it is easy to check that
$S_\l\in\Phi(U_q(\nn))$.
Indeed, due to the definition of a standard tableau, no factor
of a word $w[T]$ can be of the form 
\begin{eqnarray*}
&&w[i,i,i+1], w[i,i+1,i], w[i+1,i,i], \qquad\qquad\qquad (2\le i \le r-1) \\
&&w[i,i,i-1], w[i,i-1,i], w[i-1,i,i], \qquad\qquad\qquad (2\le i \le r) \\
&&w[1,1,1,2], w[1,1,2,1], w[1,2,1,1], w[2,1,1,1],
\end{eqnarray*}
hence the only relations to check are those involving 
subwords of the type $w[i,j]$ with $|i-j|\ge 2$, that is $a_{ij}=0$.
These relations are trivially satisfied, since they 
correspond to the exchange in a standard tableau $T$
of two consecutive integers located in two cells which
are neither in the same row nor in the same column.

Let us prove that $S_\l=\Phi(b^*_\l)$ by induction
on $|\l|=\l_1+\cdots +\l_k$.
This holds trivially for $|\l|=0$.
Now, it is clear from the definition of $\ep_i$ that we
have
\[
\ep_i(S_\l) = 
\left\{
\matrix{
S_{(\l_1,\ldots ,\l_j-1,\ldots ,\l_k)} &\mbox{\rm if $\l_j= i$ and
  $\l_{j+1}\not = i-1$,}\cr
0\hfill &\mbox{otherwise.}\hfill\hfill
}
\right.
\] 
Therefore, by induction, 
$\ep_i(S_\l)=\Phi(e'_i(b^*_\l))$ 
for all $i=1,\ldots ,r$.
Hence $\Phi^{-1}(S_\l)=b^*_\l$ by (\ref{eq:E'}).
Finally, since $\max (S_\l) = g(\l)$ we have $S_\l=b^*_{g(\l)}$.
\cqfd

More generally, we can consider skew shifted Young diagrams $\l/\m$,
where $\l=(\l_1>\cdots >\l_j>0)$ and $\m=(\m_1>\cdots >\m_k>0)$
are strict partitions with $j\ge k$ and $\l_i\ge\m_i$ for 
$i=1,\ldots,k$.
Then as above we define shifted standard Young tableaux $T$ of 
shape $\l/\m$ and we associate to them a word $w(T)$ obtained 
by reading the contents of $\l/\m$ in the order specified by~$T$.
Finally we set
$
g(\l/\m) = w[\m_1+1,\ldots,\l_1,\m_2+1,\ldots,\l_2,\ldots
,\m_j+1,\ldots,\l_j].
$
Here we understand that $\m_i+1,\ldots,\l_i$ is the empty string
if $\m_i=\l_i$, and $\m_i=0$ for $i>k$. 

The next proposition generalizes Proposition~\ref{chartabshif} to skew
shifted Young diagrams.
\begin{proposition}\label{chartabshif2}
We have 
\begin{equation}\label{tabl2}
b^*_{g(\l/\m)} = \sum_T w[T]\,,
\end{equation}
where $T$ runs through the set of all standard shifted Young
tableaux of shape $\l/\m$.
\end{proposition}
\proof
Recall the anti-automorphism $\tau$ of $\F$ defined in Proposition~\ref{proptau}.
It induces on $U_q(\nn)$ the anti-automorphism (also denoted by
$\tau$) which fixes the Chevalley generators $e_i$.
We shall use the following known properties of $\Bs$:

(a) $\tau(\Bs) = \Bs$; 

(b) given $b^*\in\Bs$ if we have 
$(e'_i)^k(b^*)\not = 0$ and $(e'_i)^{k+1}(b^*) = 0$, then
$(e'_i)^{(k)}(b^*)\in\Bs$.

\noindent
Let ${e'}^\dag_i=\tau\circ e'_i\circ\tau$.
Combining (a) and (b) we get:

(c) given $b^*\in\Bs$ if we have 
$({e'}^\dag_i)^k(b^*)\not = 0$ and $({e'}^\dag_i)^{k+1}(b^*) = 0$, then
$({e'}^\dag_i)^{(k)}(b^*)\in\Bs$.

\noindent
Let us argue by induction on $|\m|=\sum_i\m_i$.
If $|\m|=0$, the result is true by Proposition~\ref{chartabshif}.
Suppose now that the result holds for all $\l/\nu$ with
$|\nu|=p$, and choose $\m$ with $|\m|=p+1$.
There exists a strict partition $\nu$  
with $|\nu|=p$ contained in $\l$ such that we pass from $\l/\nu$
to $\l/\m$ by erasing one cell situated at the left end of its row.
Let $i$ be the content of this cell.
Recall from the proof of Lemma~\ref{factor}
that in the shuffle realization ${e'}^\dag_i$ act as $\ep^\dag_i$,
that is by removing the first letter if it is equal to $w[i]$
and by zero otherwise.
It is then easy to check that $(\ep^\dag_i)^k$ applied to
$b^*_{g(\l/\nu)}$ is zero for $k>1$
and is equal to the right-hand side of (\ref{tabl2}) for $k=1$.
Thus the statement follows from~(c).
\cqfd 

\subsection{} 
In this section we propose a conjectural type $B$ analogue of Theorem~\ref{qcharA}.

\subsubsection{}
Let ${\cal H}_m(t)$ denote the affine Hecke-Clifford superalgebra
defined by Jones and Nazarov \cite{JN} and further studied by
Brundan and Kleshchev \cite{BK1}.
We assume that the ground field is $\C$ and that the quantum
parameter $t\in\C^*$ is not a root of $1$.

We shall not write the full presentation of ${\cal H}_m(t)$,
but only recall that it consists of even generators
$T_1,\ldots , T_{m-1},X_1,\ldots ,X_m$
together with odd generators $C_1,\ldots ,C_m$, and that the 
$X_i$ are invertible and pairwise commutative. 
Brundan and Kleshchev have introduced a class of finite-dimensional
${\cal H}_m(t)$-modules, called integral.
These are the modules on which all eigenvalues of
$X_1+X_1^{-1},\ldots , X_m+X_m^{-1}$ are of the form
\[
t(i)=2\,{t^{2i-1}+t^{-2i+1}\over t+t^{-1}}\,, \qquad (i\in\N^*).
\]
Fix $r\ge 2$ and 
let $\CC_{m,r}$ denote the category of integral 
${\cal H}_m(t)$-modules for which these eigenvalues
belong to the finite subset $\{t(1),\ldots , t(r)\}$.
Let $\RR = \bigoplus_{m\in\N} \RR_{m,r}$, where $\RR_{m,r}$
is the complexified Grothendieck group of $\CC_{m,r}$.
As in \ref{charshuffle}, $\RR$ is endowed with a  
multiplication $\times$ coming from a modification $\circledast$
of parabolic induction appropriate to the superalgebra setting
\cite{BK1}.

There are $r$ irreducible modules $L(1), \ldots , L(r)$
in $\CC_{1,r}$, and they are all of dimension $2$. 

\subsubsection{}
Recall the discussion of \ref{special1}.
Let $\xi_1,\ldots ,\xi_r$ be the elements of $\C[N]$
obtained by specializing at $q=1$ the Chevalley generators
$e_1,\ldots , e_r$.
Brundan and Kleshchev have proved that
there exists an algebra isomorphism
from $\C[N]$ to $\RR$ which maps $\xi_i$ to $[L(i)]$.
Moreover there is a natural labelling of the basis of $\RR$ 
consisting of the classes of simple modules by the vertices
of the crystal graph of $U_q(\nn)$.
Here, `natural' means that the Kashiwara operators on the crystal
correspond to taking the socle of the $i$-restriction of a simple
module.

Brundan and Kleshchev have also introduced a notion of character
for the integral modules. 
Let $M$ be a module in $\CC_{m,r}$, and let $M[i_1,\ldots ,i_m]$
denote the generalized eigenspace of the pairwise commuting operators
$X_1+X_1^{-1},\ldots, X_m+X_m^{-1}$ corresponding to the 
eigenvalues $t(i_1),\ldots ,t(i_m)$, respectively. 
Write $\delta(i_1,\ldots , i_m)$ for the number of occurences
of $1$ in the list $(i_1,\ldots , i_m)$.
Then
\[
\ch M = \sum_{1\le i_1,\ldots ,i_m\le r} 2^{\lfloor\delta(i_1,\ldots ,
  i_m)/2\rfloor-m}
\,\dim M[i_1,\ldots ,i_m]\,w[i_1,\ldots ,i_m]\,.
\] 
This is an element of $\F_\C$. 
Moreover, as in \ref{charshuffle}, there holds \cite{BK1}
\[
\ch( M_1\circledast M_2) = \ch M_1 \shuffle \ch M_2\,.
\]
\begin{conjecture}
Let $b^*\in\Bs$ be an element of principal degree $m$. 
The specialization at $q=1$ of
$\Phi(b^*)$ is the character of an irreducible integral
${\cal H}_m(t)$-module.
\end{conjecture} 
This conjecture is supported by the calculations of 
Lemma~\ref{examplesB}, which agree with the character calculations
of Brundan and Kleshchev (\cite{BK1}, 5-f), and
by Proposition~\ref{chartabshif}, which agrees with the 
known characters of the finite Hecke-Clifford superalgebras
introduced by Olshanski \cite{O}.
We believe that the $b^*_{g(\l/\m)}$ for $|\l/\m|=m$
give the complete list of irreducible 
integral `tame' characters of ${\cal H}_m(t)$, \ie
the characters of the integral simple modules on which 
$X_1+X_1^{-1},\ldots, X_m+X_m^{-1}$
act semi-simply.
Finally, (\ref{GLB})
suggests that the representations of the affine Hecke-Clifford
superalgebras corresponding to the vectors
$
b^*_{w[i,\ldots,j]}\ (1\le i \le j)$ and
$b^*_{w[1,\ldots,j,1,\ldots,k]}\ (1\le j < k)
$
should play the role of the `segment' representations
in the Zelevinsky classification of irreducible representations
of affine Hecke algebras.


\section{Good Lyndon words and root vectors} \label{SECT6}

We give below the description of the root vectors 
$b^*_l=E^*_l$ for all root systems 
except $F_4$ and $G_2$, for the standard total ordering of $I$, that is, 
$w[1]<w[2]<\cdots <w[r]$.
(For type $G_2$, see \ref{EXG2}.)
For types $A, B, C, D$ we provide a closed
$q$-shuffle formula for the root vectors,
and for types $A, D, E$
we give a simple combinatorial formula
(Proposition~\ref{coeff1}).

For $l=w[i_1,\ldots ,i_k]\in\GL$, we write 
$b^*[i_1,\ldots ,i_k]$ rather than $b^*_l$.
\subsection{} {\it Type $A_r$}. \label{TypeA}
The simple roots are numbered as shown on the following Dynkin diagram:
\begin{center}
\leavevmode
\epsfxsize =6.8cm
\epsffile{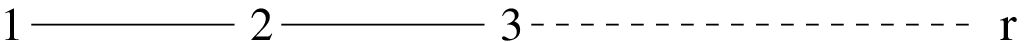}
\end{center}
The set of good Lyndon words is
$
\GL=\{w[i,i+1,\ldots , j],\ 1\le i\le j\le r\},
$
and the corresponding root vectors are 
\[
b^*[i,i+1,\ldots ,j] = w[i,i+1,\ldots ,j],
\qquad
(1\le i\le j\le r),
\]
as can be checked easily by induction on $j-i$, using formula 
(\ref{kappal}) for $b^*_l = E^*_l$.
\subsection{} {\it Type $B_r$}.\label{TypeB} 
We choose the numbering
\begin{center}
\leavevmode
\epsfxsize =6.8cm
\epsffile{dynkinBn.eps}
\end{center}
of the Dynkin diagram. 
The set of good Lyndon words is
\[
\GL=\{w[i,\ldots,j],\ 1\le i\le j\le r\}
\cup
\{w[1,\ldots,j,1,\ldots,k],\ 1\le j<k\le r\}.
\]
As in \ref{TypeA}, we have
\begin{equation}\label{rootfacile}
b^*[i,\ldots,j] = w[i,\ldots,j],
\qquad
(1\le i\le j\le r).
\end{equation}
The other root vectors are given by
\begin{lemma} \label{formularoot}
For $1\le j<k\le r$, one has
\[
b^*[1,\ldots , j, 1, \ldots , k] =
[2]_1\,w[1]\cdot(w[2,\ldots ,j]\shuf w[1,\ldots ,k]),
\]
where $\cdot$ denotes the concatenation product.
If $j=1$ we understand $w[2,\ldots ,j] = w[\,]$.
\end{lemma}
\proof
By (\ref{rootfacile}), $w[1,\ldots ,k]$ and $w[2,\ldots ,j]$
belong to $\U$. It follows that $w[2,\ldots ,j]\shuf w[1,\ldots ,k]$
also belongs to $\U$.
By Theorem~\ref{th2}, we can see now that
$f=w[1]\cdot(w[2,\ldots ,j]\shuf w[1,\ldots ,k])$ belongs to $\U$.
Indeed, we only have to check those equations 
(\ref{eqembed}) involving the first letter $w[1]$ of all words
occuring in $f$, that is, those equations for which $z=w[\,]$,
$i=1$ and $j=2$.
Since there are only $2$ occurences of $w[1]$ in each word  
and $2<1-a_{12}=3$, there are in fact no new relations
to check.
It is easy to see that $\max(f) = w[1,\ldots,j,1,\ldots,k]$,
hence, by Proposition~\ref{smallest},  $f$ is
proportional to $b^*[1,\ldots,j,1,\ldots,k]$. 
Finally, we have to show that the
proportionality factor $\gamma$ is equal to $[2]_1$.
Write 
$
l=w[1,\ldots,j,1,\ldots,k]
$
and let $l=l_1l_2$ be the co-standard factorization.
If, $k=j+1$ then $l_1=w[1,\ldots,j]$ and $l_2=w[1,\ldots,k]$.
Combining (\ref{kappal}) with (\ref{rootfacile}), we can
calculate $\gamma=[2]_1$.
If $k>j+1$ then $l_1=w[1,\ldots,j,1,\ldots ,k-1]$ and $l_2=w[k]$,
so we can show by induction on $k-j$ that $\gamma=[2]_1$.
\cqfd

\subsection{} {\it Type $C_r$}.
We choose the numbering
\begin{center}
\leavevmode
\epsfxsize =6.8cm
\epsffile{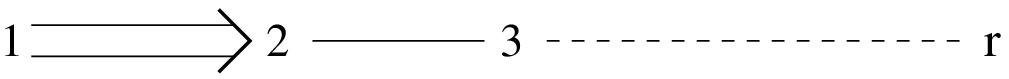}
\end{center}
The set of good Lyndon words is
\[
\GL=\{w[i,\ldots,j],\ 1\le i\le j\le r\}
\cup
\{w[1,\ldots,k,2,\ldots,j],\ 1<j\le k\le r\}.
\]
As in \ref{TypeA}, we have
$
b^*[i,\ldots,j] = w[i,\ldots,j]
\
(1\le i\le j\le r).
$
The other root vectors are given by the following lemma,
whose proof is similar to that of Lemma~\ref{formularoot} and will be
omitted.
\begin{lemma}
For $2\le j\le k\le r$, one has 
$
b^*[1,\ldots ,k,2,\ldots ,j] = w[1]\cdot(w[2,\ldots ,j]\shuf w[2,\ldots ,k]). 
$
\end{lemma}

\subsection{} {\it Type $D_r$}.\label{TypeD}
We choose the numbering
\begin{center}
\leavevmode
\epsfxsize =6.8cm
\epsffile{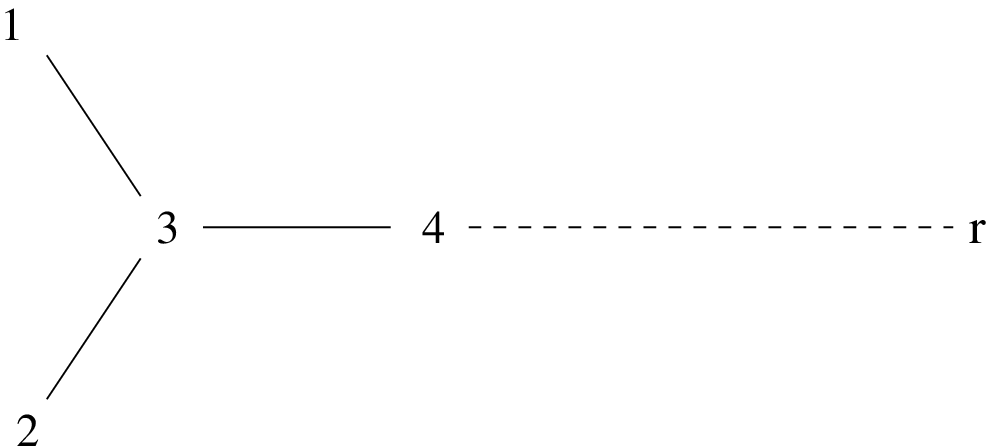}
\end{center}
The set of good Lyndon words is
\[
\GL=
\{w[1]\} \cup
\{w[1,3,\ldots,i],\ 3\le i\le r\}
\ \cup 
\qquad\qquad\qquad\qquad\qquad
\]
\[
\qquad\qquad\qquad\qquad
\{w[i,\ldots,j],\ 2\le i\le j\le r\}
\cup
\{w[1,3,\ldots,k,2,\ldots,j],\ 2\le j < k\le r\}.
\]
As in \ref{TypeA}, we have 
\[
b^*[1,3,\ldots,i] = w[1,3,\ldots,i],\ \ (3\le i\le r),\qquad
b^*[i,\ldots,j] = w[i,\ldots,j],
\ \
(2\le i \le j\le r).
\]   
The other root vectors are given by
\begin{lemma}
For $2\le j< k\le r$,
\[
b^*[1,3,\ldots,k,2,\ldots,j] = 
w[1]\cdot(w[2,\ldots ,j]\shuf w[3,\ldots ,k] - q\,w[2,\ldots ,k]\shuf
w[3,\ldots ,j]),
\]
where for $j=2$, we understand $w[3,\ldots , j] = w[\,]$.
\end{lemma}
\proof
The proof is similar to that of Lemma~\ref{formularoot}.
First, we see as above that 
\[
u=w[2,\ldots ,j]\shuf w[3,\ldots ,k] - q\,w[2,\ldots ,k]\shuf
w[3,\ldots ,j]
\]
belongs to $\U$. 
Secondly, one can check that all words occuring 
in $u$ start with the letter $w[3]$, since all words starting
with $w[2]$ cancel out.
Moreover, no word starts with $w[3,3]$.
Therefore, $f=w[1]\cdot u$ also belongs to $\U$ (no new
relations to be checked) and the proof is concluded as above.
\cqfd

\subsection{} {\it Type $E$}.\label{TypeE}
We choose the same numbering as in \cite{LR}.
The set $\GL$ can then be read from the Lyndon paths in the $E_8$-tree
of \cite{LR}.
For example the good Lyndon word associated to the highest root of 
$E_8$ is 
\[
w[1,3,4,5,6,7,8,2,4,5,6,3,4,5,2,4,3,1,3,4,5,6,7,8,2,4,5,6,7].
\]
\subsection{}
For types $A, D, E$ we have the following combinatorial description
of the root vectors attached to the good Lyndon words above.

Let $\sim$ be the equivalence relation in $\M$ defined
by $w\sim w'$ if only if $w'$ can be obtained from 
$w$ by a sequence of commutations of two adjacent letters
$w[i]$ and $w[j]$ with $a_{ij}=0$.

\begin{proposition}\label{coeff1}
Let $\g$ be of type $A_r, D_r$ or $E_r$.
For any $l\in\GL$, we have
$
b^*_l = \sum_{w\sim l} w\,.
$
In particular, all words occuring in $b^*_l$ have coefficient $1$,
and $\kappa_l=1$.
\end{proposition}
\proof
We see by inspection of the sets $\GL$ given above for types
$A,D,E$ that no $l\in\GL$ has a factor of the form
\begin{equation}\label{iij}
w[i,i,j], \quad w[i,j,i], \quad w[j,i,i]
\end{equation}
with $a_{ij}=-1$.
Moreover, for any $i$ occuring more than once in $l$, and any
factor $x=w[i]\cdot y\cdot w[i]$ of $l$, we can check that $y$ contains
at least $2$ letters $w[j]$ and $w[k]$ with $a_{ij}=a_{ik}=-1$.
This implies that no $w$ equivalent to $l$ has a factor
of the form (\ref{iij}). 
It follows that $f_l= \sum_{w\sim l} w$
satisfies the equations of Theorem~\ref{th2}, and therefore 
belongs to $\U$. 
Moreover, again by Theorem~\ref{th2}, any $u\in\U$
in which the word $l$ occurs with coefficient $\gamma$
contains all words of $f_l$ with the same coefficient $\gamma$.
Hence, $f_l$ is equal to the element $m_l$ of the
basis $\{m_g\}$ introduced in the proof of Proposition~\ref{monbasis},
and by Proposition~\ref{smallest}, $f_l$ is proportional to $b^*_l$. 

It remains to prove that the coefficient of proportionality is equal
to $1$.
By (\ref{kappal}), this amounts to prove that the coefficient 
of $l$ in 
\[
{(-1)^{\ell(l)-1}\over q^{N(|l|)}(E_l,E_l)}\,r_l
=
{1\over (q-q^{-1})^{\ell(l)-1}}\,r_l
\]
is equal to $1$ (in the simply laced case, we have
$(E_l,E_l)=(1-q^2)^{\ell(l)-1}$ and $N(|l|)=1-\ell(l)$). 
To see this we proceed by induction on $\ell(l)$ and
consider the co-standard factorization
$l=l_1l_2$ of $l$.
Since 
$r_l=r_{l_1}\shuf r_{l_2} - q^{(|l_1|,|l_2|)} r_{l_2}\shuf r_{l_1}$,
we are reduced to prove that the coefficient of
$l$ in $l_1\shuf l_2 - q^{(|l_1|,|l_2|)} l_2\shuf l_1$ is equal
to $q-q^{-1}$. 
To show this, it is enough by Proposition~\ref{propo2} to 
show that the coefficient of $l$ in $l_1\shuf l_2$ is equal to $q$.
By (\ref{eq:1.5}) (\ref{eq:1.6}) this coefficient is equal to 
$q^{-(|l_1|,|l_2|)}$, so all we have to prove is that 
$(|l_1|,|l_2|)=-1$.

For types $A$ and $D$ we see immediately from \ref{TypeA} and \ref{TypeD}
that $l_2$ is always reduced to the last letter of $l$,
and we can easily check that $(|l_1|,|l_2|)=-1$.
For example in type $D_r$, if $l=w[1,3,\ldots,k,2,\ldots ,j]$
with $2\le j<k \le r$, we have 
\[
(|l_1|,|l_2|)=(\a_{j-1}+\a_j+\a_{j+1}+\a_{j-1}\,,\,\a_j)= 
-1+2-1-1=-1\,.
\]  
For type $E$, the equality $(|l_1|,|l_2|)=-1$ can be checked 
from the lists of good Lyndon words given in \cite{LR}.
(In most of the cases $l_2$ is reduced to the last letter of $l$
and the calculation of $(|l_1|,|l_2|)$ is very easy.)
\cqfd

We believe that Proposition~\ref{coeff1} also holds
for all other total orderings of the set of simple roots.

\bigskip
\centerline{\bf Acknowledgements}

\medskip\noindent
I thank A. Berenstein, E. Frenkel, M. Rosso and A. Zelevinsky
for helpful discussions.
This work was done during a visit at the National Technical
University of Athens.
I am grateful to the Mathematics Department of N. T. U. A. for
its hospitality.
  

%
\bigskip
\small

\noindent
\begin{tabular}{ll}
{\sc B. Leclerc} : &
D\'epartement de Math\'ematiques,
Universit\'e de Caen, Campus II,\\
& Bld Mar\'echal Juin,
BP 5186, 14032 Caen cedex, France\\
&email : {\tt leclerc@math.unicaen.fr}
\end{tabular}
\end{document}